\pdfoutput=1
\documentclass{amsart}
\usepackage{amsmath,amssymb}
\usepackage{epsfig}
\usepackage{amsthm}
\usepackage{url}
\usepackage{enumerate}

\newcommand{\N}{\mathbb{N}}

\newcommand{\R}{\mathbb{R}}
\newcommand{\C}{\mathbb{C}}

\newtheorem{theorem}{Theorem}[section]

\newtheorem{remark}[theorem]{Remark}

\newtheorem{algorithm}{Algorithm}

\usepackage{color}
\definecolor{darkblue}{rgb}{0,0,0.6}
\definecolor{darkgreen}{rgb}{0,0.6,0}
\definecolor{darkred}{rgb}{0.6,0,0}
\definecolor{orange}{rgb}{1,0.7,0}
\definecolor{darkorange}{rgb}{1,0.4,0.3}
\definecolor{brown}{rgb}{0.7,0.3,0}
\definecolor{MatBlue}{rgb}{0,0.4431,0.7373}
\definecolor{MatOrange}{rgb}{0.9255,0.6902,0.1216}
\definecolor{MatRed}{rgb}{0.8471,0.3216,0.0941}
\definecolor{revred}{rgb}{0.9,0,0.7}

\newcommand{\domain}{\Omega}
\newcommand{\modeldomain}{\Omega_0}

\newcommand{\sigbar}{\sigma}

\newcommand{\sig}{\tilde{\sigma}}
\newcommand{\sigmean}{\tilde{\sigma}_{\ast}}

\newcommand{\sigbarsamp}[1]{\sigma^{(#1)}}
\newcommand{\alpbarsamp}[1]{\vartheta^{(#1)}}

\newcommand{\Vmeas}{\mathcal{V}}

\newcommand{\Uforw}{\mathcal{U}}
\newcommand{\Uacc}{\Uforw (\sigbar,z,\vartheta)}
\newcommand{\Umod}{\Uforw_0 (\sig,z)}
\newcommand{\Uaccsamp}[1]{\Uforw (\sigbarsamp{#1},z^{(#1)},\alpbarsamp{#1})}
\newcommand{\Umodsamp}[1]{\Uforw_0 (\sigmean,z^{(#1)})}

\newcommand{\etamean}{\eta_{\ast \vert \omega}}
\newcommand{\etacov}{\Gamma_{\eta \vert \omega}}

\newcommand{\epsmean}{\varepsilon_{\ast}}
\newcommand{\epscov}{\Gamma_\varepsilon}
\newcommand{\epssamp}[1]{\varepsilon^{(#1)}}
\newcommand{\transp}{^{\rm T}}

\begin{document}
\title[Approximation error method for head imaging by EIT]{Approximation error method for imaging the human head by electrical impedance tomography}

\author{V.~Candiani}
\address{Aalto University, Department of Mathematics and Systems Analysis, P.O.~Box 11100, FI-00076 Aalto, Finland} 
\email{valentina.candiani@aalto.fi}

\author{N.~Hyv\"onen}
\address{Aalto University, Department of Mathematics and Systems Analysis, P.O. Box 11100, FI-00076 Aalto, Finland} 
\email{nuutti.hyvonen@aalto.fi}

\author{J.~P.~Kaipio}
\address{University of Auckland, Department of Mathematics, Private Bag 92019, Auckland 1142, New Zealand} 
\email{jari@math.auckland.ac.nz}

\author{V.~Kolehmainen}
\address{University of Eastern Finland, Department of Applied Physics, Kuopio Campus, P.O.~Box 1627, FI-70211 Kuopio, Finland} 
\email{ville.kolehmainen@uef.fi}

\thanks{This work was supported by the Academy of Finland and the Jane and Aatos Erkko Foundation.}



\begin{abstract}
 This work considers electrical impedance tomography 
 imaging of the human head, with the ultimate goal of
 locating and classifying a stroke in emergency care. One of the main difficulties in the envisioned application is that the electrode locations and the shape of the head are not precisely known, leading to significant imaging artifacts due to impedance tomography being sensitive to modeling errors. In this study, the natural variations in the geometry of the head and skull are modeled based on a library of head anatomies. The effect of these variations, as well as that of misplaced electrodes, on (absolute) impedance tomography measurements is in turn modeled by the approximation error method. This enables reliably reconstructing the conductivity perturbation caused by the stroke in an average head model, instead of the actual head, relative to its average conductivity levels. The functionality of a certain edge-preferring reconstruction algorithm for locating the stroke is demonstrated via numerical experiments based on simulated three-dimensional data. 
\end{abstract}

\maketitle

\section{Introduction}
\label{sec:intro}

{\em Electrical impedance tomography} (EIT) employs a set of contact electrodes to drive alternating electric currents into an examined physical body and to measure the resulting electrode potentials. The aim is to reconstruct (useful information on) the conductivity distribution inside the body of interest, which corresponds to a nonlinear and highly ill-posed inverse elliptic boundary value problem; see~\cite{Borcea02,Cheney99,Uhlmann09} for more information. This work considers stroke detection and classification between ischemic and hemorrhagic strokes by EIT. The leading idea is to combine a principal component model for the variations in the anatomy of the human head \cite{Candiani19,Candiani20} with the approximation error approach \cite{Kaipio05,kaipio2013} and a certain edge-promoting reconstruction algorithm \cite{Harhanen15} in order to introduce an image reconstruction method that is robust with respect to geometric uncertainties in the measurement configuration and carries potential to identify stroke events despite the partial shielding by the skull.

EIT is known to be highly sensitive to mismodeling of the measurement setup~\cite{Barber88,Breckon88,Kolehmainen97}, which is a major hindrance for many medical applications where the exact shape of the imaged object and the precise positions of the electrodes are often unknown. In particular, this applies to stroke detection and classification in emergency care, where EIT data has to be collected and reconstructed without access to an auxiliary scan of the patient's head and the electrode geometry.
The simplest (partial) cure for inaccurate geometric modeling is difference imaging~\cite{Barber84}: if electrode measurements are taken at two separate times (or at two frequencies), the modeling errors cancel out (partially) when the change in the conductivity is reconstructed based on the difference of the two measurements and a linearized approximation of the forward model. However, the linearization approach leads to suboptimal results, and moreover, reference images are not always available. Particularly, in stroke classification one does not have access to a reference measurement taken before the occurrence of the stroke to be classified.

Previous techniques for handling geometric uncertainties in the fully nonlinear setting of EIT are listed in the following. The first one was introduced in \cite{Kolehmainen05,Kolehmainen07}, where the mismodeling of the measurement geometry was compensated by reconstructing a suitable, mildly anisotropic conductivity. The papers \cite{Darde12,Darde13a,Darde13b,Hyvonen17c} adopted a conceptually straightforward approach to coping with incomplete information on the measurement setup: the derivatives of the electrode potentials with respect to the body shape and the electrode positions were utilized in a regularized Newton-type algorithm to simultaneously reconstruct all (geometric) unknowns in the considered setting. This method was extended to (three-dimensional) head imaging in \cite{Candiani19}, but without properly modeling the resistive skull layer. Polynomial surrogates for the dependence of the electrode measurements on unknown geometric parameters were built and employed in straightforward and fast Tikhonov regularization for two-dimensional EIT in~\cite{Hyvonen17}.

In this work we combine the principal component head model from \cite{Candiani19,Candiani20} with a `hybrid' reconstruction approach, taking elements from both the Bayesian framework, i.e.~approximation errors and likelihood models, and the regularization framework, i.e.~an edge-promoting penalty term with the level of regularization tuned according to the Morozov principle. The key idea of the approximation error method~\cite{Kaipio05}, which has previously been applied to EIT with geometric uncertainties in \cite{Nissinen11,Nissinen11b,kaipio2013,Hadwin14},
is to model the error caused by the inaccurate information on the electrode positions and the anatomy of the examined patient head as an additive approximation error term in the measurement model, and then to marginalize the likelihood model of the measurements with respect to the approximation error to treat  the model uncertainties. A computationally plausible closed form solution for the marginalization can be obtained by approximating the joint distribution of the approximation error and the conductivity by a Gaussian model. 

The (Gaussian) second order statistics of the approximation error are estimated by a simple machine learning procedure. The learning data is formed by computing simulated EIT measurements for an ensemble of random head and skull shapes and electrode locations. The sought-for statistics are then obtained by considering the deviation of these measurements from the ones corresponding to the intended electrode positions and the mean anatomy inside the average head of the employed principal component model~\cite{Candiani20}. The actual inversion is performed by reconstructing the deviation of the conductivity (inside the examined head) from an expected conductivity
{\em using the average head model}\footnote{The considered inverse problem indeed corresponds to absolute imaging, as we do not have a reference voltage data for the expected state --- the conductivity is just parametrized by its deviation from the prior mean.}, accounting for both the actual measurement noise and the approximation error noise within the Bayesian paradigm.

The prior information that the stroke presumably only affects a certain limited region inside the brain is incorporated by introducing an edge-preferring prior density for the (discretized) change in the conductivity compared to the mean anatomy. Computing a {\em maximum a posteriori} (MAP) estimate, i.e.~the maximizer of the posterior density, thus corresponds to finding a minimizer for a nonlinear Tikhonov-type functional that is nonquadratic in both the discrepancy and the penalty term. In the terminology of regularization theory, this corresponds to applying Tikhonov regularization to a nonlinear inverse problem with a suitable nonlinear penalty term, such as {\em total variation} (TV)~\cite{Rudin92} or {\em Perona--Malik} (PM)~\cite{Perona90}.

We choose a slightly modified version of the algorithm introduced in \cite{Harhanen15} for tackling the aforementioned nonquadratic Tikhonov functional; see~\cite{Arridge13} for the original ideas behind this technique for linear inverse problems as well as \cite{Hannukainen15,Hannukainen16} for the application of (essentially) the same algorithm to diffuse optical tomography and quantitative photoacoustic tomography. The method consists of two nested iterations: The outer loop corresponds to sequential linearizations of the forward map.  Although the inner iteration is motivated by taking a couple of lagged diffusivity steps~\cite{Vogel96} for a linear inverse problem with an edge-promoting penalty term, the resulting quadratic Tikhonov functionals are at the end replaced by priorconditioning~\cite{Calvetti07,Calvetti12,Calvetti05,Calvetti08} and a subsequent application of a certain Krylov subspace iteration, namely LSQR~\cite{Paige82a,Paige82b}. An early stopping rule serves as the only source of regularization for the LSQR iterations, with the approximation error noise taken appropriately into account when implementing the Morozov discrepancy principle. Other (edge-promoting) algorithms could presumably be used in place of that from~\cite{Harhanen15} in our numerical experiments; the main reason for relying on the ideas in \cite{Harhanen15} is the capability of the associate algorithm to produce good quality reconstructions on dense and unstructured {\em finite element} (FE) meshes in a matter of minutes on a standard laptop computer.

This text is organized as follows. Section~\ref{sec:forward} recalls the {\em complete electrode model} (CEM) \cite{Cheng89,Somersalo92} that is the most accurate model known for real-world EIT measurements. Section~\ref{sec:param} discusses the principal component parametrization for the human head, introduced originally in \cite{Candiani19} and supplemented with a treatment of the skull in \cite{Candiani20}. The approximation error methodology and the edge-promoting reconstruction algorithm are reviewed in Sections~\ref{sec:approx} and \ref{sec:algo}, respectively. The numerical experiments based on simulated data are presented in Section~\ref{sec:numer}. Finally, Section~\ref{sec:conclusion} gives the concluding remarks.

\section{Forward model}
\label{sec:forward}

In practical EIT, $M \in \N \setminus \{ 1 \}$ contact electrodes $E_1, \dots,  E_M$ are attached to the exterior surface of a bounded Lipschitz domain $\Omega \subset \R^3$. The electrodes are assumed to be disjoint, and they are identified with the open and connected surface patches they cover. We denote $E = \cup E_m$. A single measurement of EIT corresponds to driving a current pattern $I \in \R^M_\diamond$ through the electrodes and measuring the corresponding constant  potentials $U \in \R^M_\diamond$ at the electrodes. Here,
\[ 
\R^M_\diamond = \Big\{V\in\C^M\,\Big|\, \sum_{m=1}^M V_m = 0\Big\}
\]
is the mean-free subspace of $\R^M$, to which $I$ belongs due to conservation of charge and $U$ via choosing the ground level of potential appropriately. The contact resistances at the electrode-object interfaces are modeled by a vector $z  \in \R_+^M$. 

 The mathematical model that most accurately predicts real-life EIT measurements is the {\em complete electrode model} (CEM) \cite{Cheng89}: The electromagnetic potential $u$ inside $\Omega$ and the potentials on the electrodes $U$ satisfy
\begin{equation}
\label{eq:cemeqs}
\begin{array}{ll}
\displaystyle{\nabla \cdot\sigma\nabla u = 0 \qquad}  &{\rm in}\;\; \Omega, \\[6pt] 
{\displaystyle{\nu\cdot\sigma\nabla u} = 0 }\qquad &{\rm on}\;\;\partial\Omega\setminus\overline{E},\\[6pt] 
{\displaystyle u+z_m{\nu\cdot\sigma\nabla u} = U_m } \qquad &{\rm on}\;\; E_m, \quad m=1, \dots, M, \\[2pt] 
{\displaystyle \int_{E_m}\nu\cdot\sigma\nabla u\,{\rm d}S} = I_m, \qquad & m=1,\ldots,M, \\[4pt]
\end{array}
\end{equation}
interpreted in the weak sense. Here, $\nu = L^\infty(\partial \Omega, \R^n)$ is the exterior unit normal of $\partial\Omega$ and the isotropic conductivity distribution $\sigma$ is assumed to belong to
\begin{equation}
\label{eq:sigma}
L^\infty_+(\Omega) := \{ \varsigma \in L^\infty(\Omega) \ | \ {\rm ess} \inf \varsigma > 0 \}.
\end{equation}
Under the above assumptions, the potential pair $(u,U) \in H^1(\Omega) \oplus \R^M_\diamond$ is uniquely determined by \eqref{eq:cemeqs} \cite{Somersalo92}. A physical justification of \eqref{eq:cemeqs} can be found in \cite{Cheng89}.

We denote by $U(\sigma,z; I) \in \R^M_\diamond$ the functional dependence of the electrode potential pattern on the conductivity, the contact resistances and the applied electrode currents. In the case of $M-1$ linearly independent current patterns, $I_1, \dots, I_{M-1} \in \R^M_\diamond$, i.e.~a basis for $\R^M_\diamond$, we extend this notation by writing
\begin{equation}
\label{eq:forward_map}
\mathcal{U}\big(\sigma,z; I_1, \dots, I_{M-1}\big) 
= \left[U(\sigma,z; I_1)^{\rm T}, \dots, U(\sigma,z; I_{M-1})^{\rm T}\right]^{\rm T}
\in \R^{M(M-1)}.
\end{equation}
Usually the employed current basis is assumed to be known from the context, and the dependence of $\mathcal{U}$ on $I_1, \dots, I_{M-1}$ is omitted for the sake of clarity. The inverse problem of EIT consists of reconstructing (useful information about) the conductivity $\sigma$ from a noisy realization of the electrode potential measurements $\mathcal{U}(\sigma,z)$.

\section{Optimally parametrized head model}
\label{sec:param}
The three-layer head model used in this article is based on the principal component construction introduced in~\cite{Candiani19, Candiani20}. We define a {\em layer} for each considered anatomical structure: the scalp, the highly resistive skull and the interior brain. The library of $n=50$ heads from \cite{Lee16} is used for building an approximate model for the variations in the shapes and sizes of these layers over the human population. 

\subsection{Shape parametrization}
\begin{figure}
\center{
  {\includegraphics[width=5cm]{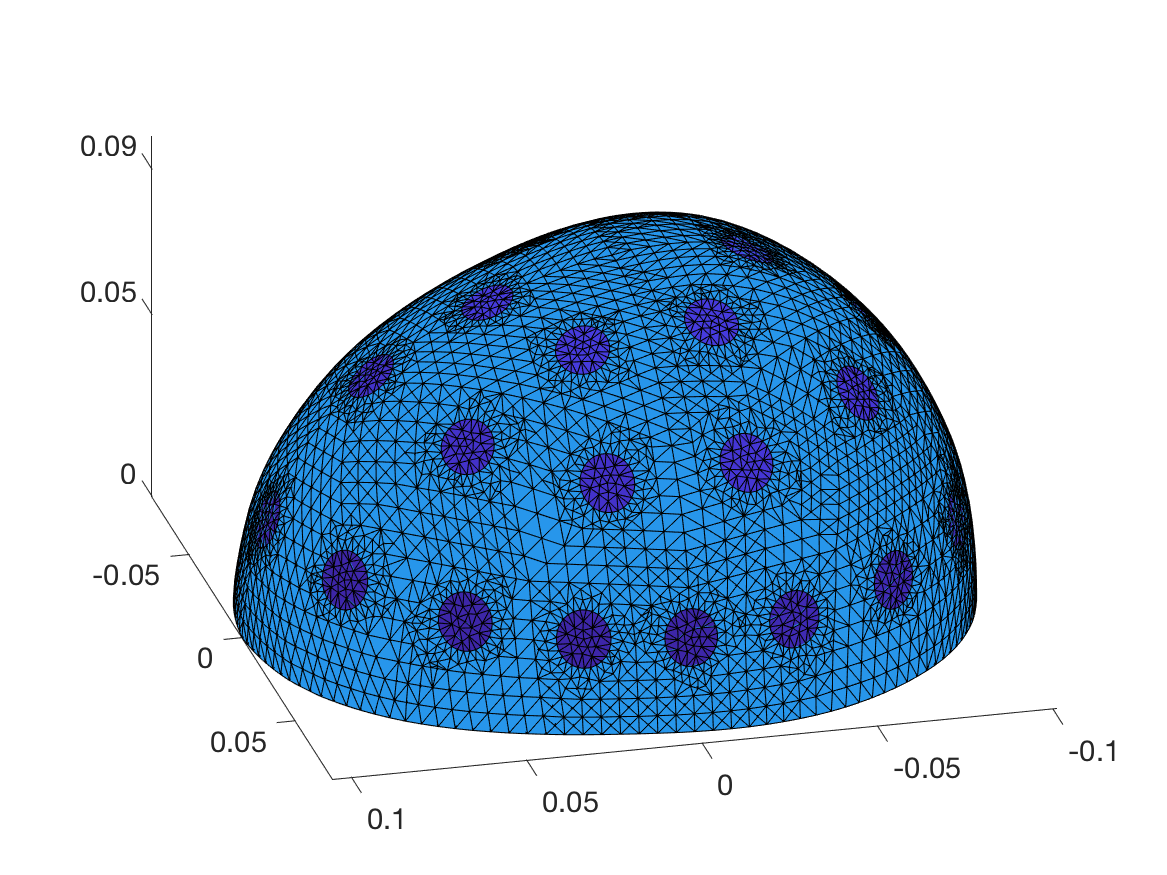}}
  \quad
  {\includegraphics[width=5cm]{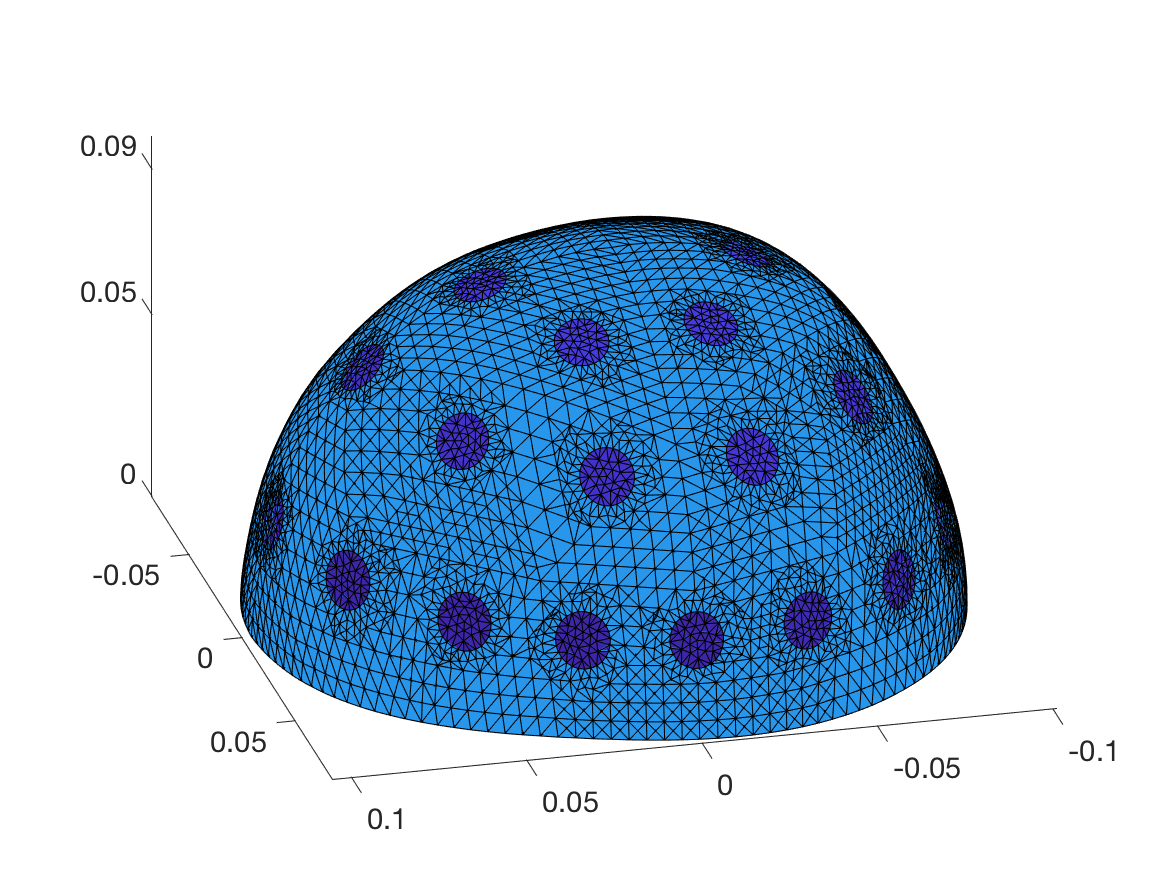}}
  }
  \center{
  {\includegraphics[width=5cm]{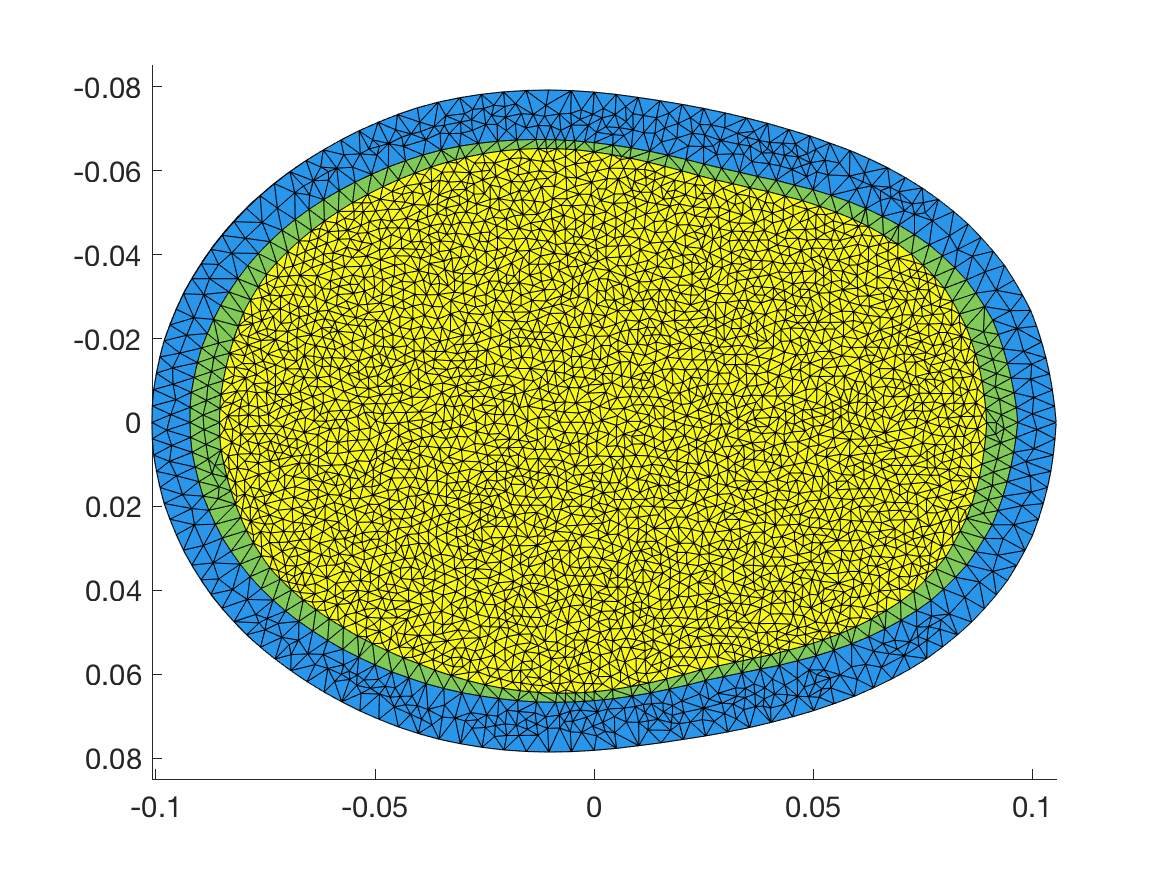}}
  \quad
  {\includegraphics[width=5cm]{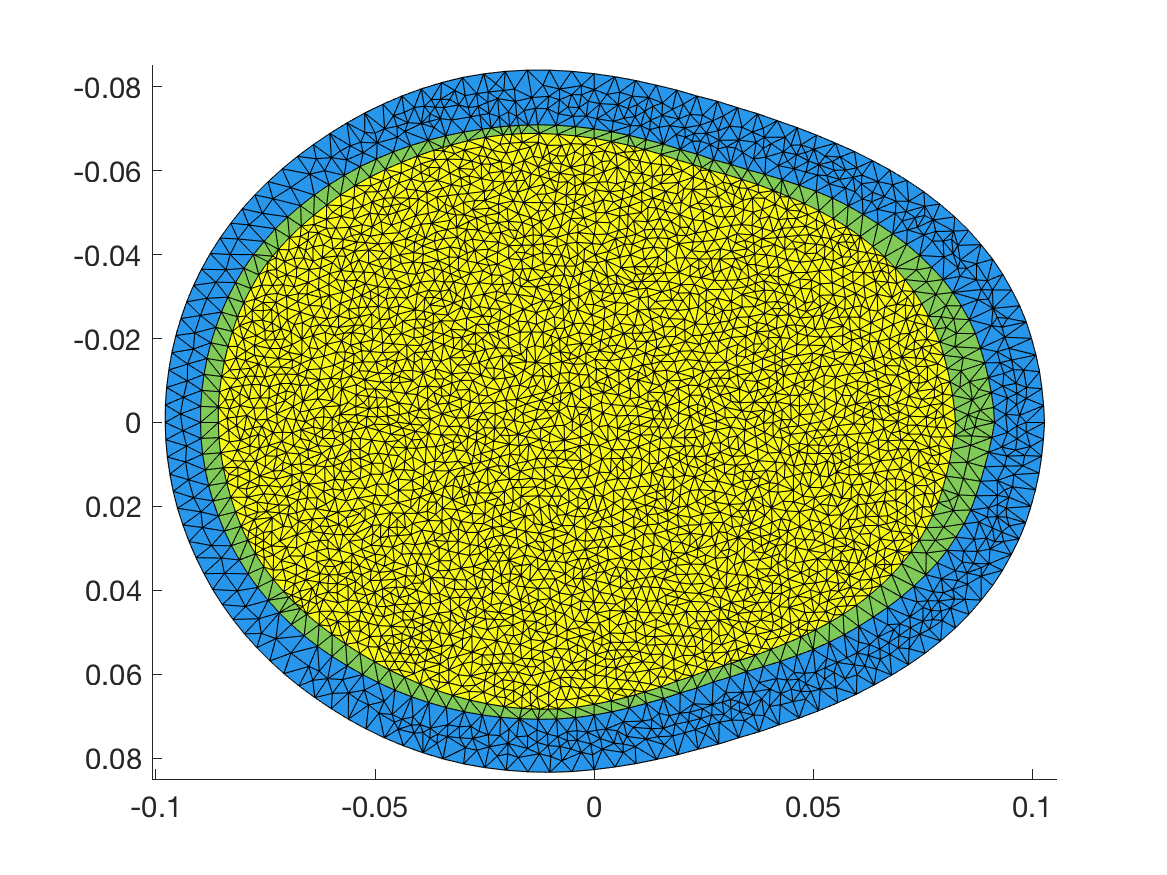}}
}
 \caption{Top row: two different head model samples with $M=32$ electrodes at their intended positions. Bottom row: the corresponding bottom faces, with the three layers associated to scalp, skull and brain tissues visible. The unit of length is meter.}
\label{fig:head}
\end{figure}

The shape of the exterior surface of the $l$th layer in the $j$th head can be represented as the graph of a function
\begin{equation}
  \label{eq:jth_head}
S_j^l:
\left\{
\begin{array}{l}
\mathbb{S}_+ \to \R^3, \\[1mm]
\hat{x} \mapsto r_j^l(\hat{x}) \, \hat{x}.
\end{array}
\right.
\end{equation}
Here $\mathbb{S}_+$ is the upper unit hemisphere, i.e.,
$$
\mathbb{S}_+ = \big\{ x \in \R^3 \; | \; | x | = 1 \ {\rm and} \ x_3 > 0 \big\},
$$
and $r_j^l:  \mathbb{S}_+ \to \R_+$ gives the distance from the origin to the exterior surface of $l$th layer as a function of the direction $\hat{x} \in \mathbb{S}_+$, with the origin set approximately at the center of mass for the bottom face of the considered $j$th head (see Figure \ref{fig:head}). Here and in what follows, we denote by $|\cdot|$ the Euclidean norm.

We define the mean for the shape of the exterior surface of the $l$th layer and the corresponding $n$ perturbations via
\begin{equation}
  \label{eq:mean_etc}
\bar{r}^l = \frac{1}{n} \sum_{j=1}^n r_j^l \qquad {\rm and} \qquad \rho_j^l = r_j^l - \bar{r}^l, \qquad j=1, \dots, n, \,\  l = 1, 2, 3.
\end{equation}
The perturbations of the three layers are represented by vector-valued linearly independent `total perturbations' $\rho_j :=[\rho_j^l]_{l=1}^3$, $j=1,\dots,n$. 

According to the above construction and assuming some regularity, a single head in the library defines an element of $[H^1(\mathbb{S}_+)]^3$, that is, a three dimensional vector whose components lie in the space $H^1(\mathbb{S}_+)$ and define the three anatomical layers. The corresponding inner product between two elements of $[H^1(\mathbb{S}_+)]^3$, say, $v$ and $w$, is defined in the natural way as
\begin{equation}
\label{eq:inner}
(v, w)_{[H^1(\mathbb{S}_+)]^3} := \sum_{l=1}^{3}(v_l, w_l)_{H^1(\mathbb{S}_+)},
\end{equation}
with $(\cdot, \cdot)_{H^1(\mathbb{S}_+)}$ denoting the standard (real) inner product of the Sobolev space $H^1(\mathbb{S}_+)$.

Our aim is to find an $\tilde{n}$-dimensional subspace $V_{\tilde{n}} \subset  [H^1(\mathbb{S}_+)]^3$, $1 \leq \tilde{n} \leq n$,
that satisfies
\begin{equation}
  \label{eq:minim}
\sum_{j=1}^{n} \min_{\eta \in V_{\tilde{n}}} \| \rho_j - \eta \|_{[H^1(\mathbb{S}_+)]^3}^2 \leq \sum_{j=1}^{n} \min_{\eta \in W} \| \rho_j - \eta \|_{[H^1(\mathbb{S}_+)]^3}^2,
\end{equation}
for all $\tilde{n}$-dimensional subspaces $W$ of $[H^1(\mathbb{S}_+)]^3$. In other words, the sought-for subspace contains on average the best approximations for the perturbations $\rho_1, \dots, \rho_n$ among the subspaces of the same dimension, with the quality of the fit measured by the squared norm of $[H^1(\mathbb{S}_+)]^3$.

Following the ideas in \cite[Lemma~3.1]{Candiani19} and \cite{Candiani20}, a subspace with the desired property can be constructed with the help of a full set of eigenpairs $\{\lambda_k, v_k\}_{k=1}^n \subset \R_+ \times \R^n$ for the symmetric positive definite matrix $R\in \R^{n \times n}$ defined componentwise as
$$
R_{ij} = (\rho_i, \rho_j)_{[H^1(\mathbb{S}_+)]^3}, \qquad i,j = 1, \dots, n.
$$
More precisely, arranging $\lambda_1 \geq \lambda_2 \geq \cdots \geq \lambda_n >0$ and choosing $v_1, \dots, v_n$ to be orthonormal,
$$
V_{\tilde{n}} = {\rm span} \{\hat{\rho}_1, \dots, \hat{\rho}_{\tilde{n}} \}, \qquad \tilde{n} \in \{1, \dots, n\},
$$
where the $[H^1(\mathbb{S}_+)]^3$-orthonormal basis function are defined via
\begin{equation}
\label{eq:rhohat}
\hat{\rho}_k := \frac{1}{\sqrt{\lambda_k}}\, [\rho_1, \dots, \rho_n] \, v_k, \qquad k=1, \dots, n,
\end{equation}
with the matrix-vector product understood as taking a linear combination of the perturbations $\rho_1, \dots, \rho_n \in [H^1(\mathbb{S}_+)]^3$. 

Our $\tilde{n}$-dimensional principal component parametrization for the exterior shape of the $l$th anatomical layer in the human head can then be written as
\begin{equation}
\label{eq:theparam_l}
S^l(\hat{x}; \alpha) = \Big( \bar{r}^l(\hat{x}) +
\sum_{k=1}^{\tilde{n}} \alpha_k \hat{\rho}_k^l(\hat{x}) \Big) \hat{x}, \qquad \hat{x} \in \mathbb{S}_+, \,\ l = 1, 2, 3,
\end{equation}
where $\hat{\rho}_k^l$ are the components of the basis functions in \eqref{eq:rhohat} and $\alpha_k \in \R$ are free shape coefficients that are common for all three layers. When considering a random three-layered head anatomy in our numerical tests, we choose to draw the vector of shape coefficients $\alpha \in \mathbb{R}^n$ from a Gaussian distribution $\mathcal{N}(0,\Gamma_\alpha)$, with the diagonal covariance matrix
\begin{equation}
  \label{eq:Gamma_alpha3}
(\Gamma_\alpha)_{pq} = \frac{1}{n-1}\sum_{j=1}^n(\rho_j, \hat{\rho}_p)_{[H^1(\mathbb{S}_+)]^3} (\rho_j, \hat{\rho}_q)_{[H^1(\mathbb{S}_+)]^3} = \frac{\lambda_p}{n-1} \delta_{pq}.
\end{equation}
See \cite{Candiani19} for further details, including the motivation behind the chosen Gaussian random model for the shape parameters.

\subsection{Generating the FEM mesh for a parametrized head}
Our workflow for generating a tetrahedral mesh for a parametrized head with a three-layer anatomical structure consists of three steps: generation of an initial surface mesh, insertion of the electrodes and tetrahedral mesh generation. The initial surface mesh is constructed by subdividing a coarse surface partition consisting of four triangles a number of times. Next, $M$ electrodes with appropriately refined meshes are inserted onto the surface with the help of the Triangle software \cite{Shewchuk96} as described in \cite{Candiani19}. The process is completed by generating a tetrahedral partition for the whole volume by TetGen \cite{Hang15} starting from the formed surface mesh and accounting for the boundaries of the anatomical layers (see Figure~\ref{fig:head}). 

\section{FEM approximation of the forward model}
\label{sec:fem}

In the rest of this work, we consider a discretized version of \eqref{eq:cemeqs}. To this end, the conductivity $\sigma$ is modeled as
\begin{equation}
\label{discr_sigma}
\sigma = \sum_{j=1}^N \sigma_j \varphi_j \, ,
\end{equation}
where $\sigma_j \in \R_+$, $j=1, \dots, N$, and $\varphi_j \in H^1(\Omega)$, $j=1, \dots, N$,  are the piecewise linear basis functions corresponding to a FE mesh of $\Omega$. We denote by $\sigma$ both the vector $\sigma \in \R^N_+$ and the corresponding conductivity defined by \eqref{discr_sigma}; we frequently use a similar identification between any vector of $\R^N$ and the corresponding function in ${\rm span} \{ \varphi_j \}_{j=1}^N \subset H^1(\Omega)$  without mentioning this explicitly, assuming the correct interpretation is clear from the context. By the solution $(u,U)$ of~\eqref{eq:cemeqs} we mean from here on the approximate FE solution in the subspace
$$
{\rm span}\{\varphi_j\}_{j=1}^N \oplus \R^M_\diamond \subset H^1(\Omega)\oplus \R^M_{\diamond};
$$
see \cite{Vauhkonen97, Vauhkonen99} for the implementation details. 

Let $\vartheta \in \R^d$ be the parameter vector defining the electrode positions as well as the three-layer anatomy of the studied head via \eqref{eq:theparam_l}; in particular, certain components of $\vartheta$ correspond to the shape coefficients $\alpha$ in \eqref{eq:theparam_l}. The forward map $\mathcal{U}: L^\infty_+(\Omega) \times \R_+^M \to \R^{M(M-1)}$, introduced originally in \eqref{eq:forward_map}, is redefined to be a finite-dimensional nonlinear map
$$
\mathcal{U}: 
\left\{
\begin{array}{l}
 (\sigma,z,\vartheta) \mapsto \mathcal{U}(\sigma,z,\vartheta), \\[1mm]
\R_+^N \times \R_+^M \times \R^d \to \R^{M(M-1)}
\end{array}
\right.
$$
in accordance with the FE discretization and the parameter vector $\vartheta$ defining the anatomical layers and the electrode configuration (for a given set of $M-1$ linearly independent current patterns).

Note that the basis functions $\varphi_j$ and their number $N$ naturally change with $\vartheta$ (and the particular FE discretization) as different domains and electrode configurations are considered, but we suppress this dependence for notational convenience. The same applies to the domain of definition for the forward map $\mathcal{U}$. An exception to this rule is the (FE discretization for the) average anatomy that corresponds to $\alpha = 0$ in \eqref{eq:theparam_l} and the electrodes at their intended positions; the corresponding total geometric parameter vector and piecewise linear FE basis functions are denoted by $\vartheta_0 \in \R^d$ and $\{\varphi_j^{(0)}\}_{j=1}^{N_0} \subset H^1(\Omega_0)$, respectively, with $\Omega_0 = \Omega(\vartheta_0)$. In particular, 
\begin{equation}
    \label{eq:inv_forv}
\mathcal{U}_0: \R^{N_0}_+ \times \R^M_+ \to \R^{M(M-1)}
\end{equation}
is our shorthand notation for the associated discretized forward map, suppressing the dependence on the parameter vector $\vartheta = \vartheta_0$.

\section{Approximation error model}
\label{sec:approx}

The approximation error approach has previously been employed for handling uncertainties related to the domain boundary in two-dimensional EIT in \cite{Nissinen11, Nissinen11b,kaipio2013}. Here, we briefly review the main idea of the approximation error model and explain how the approach is implemented for modeling errors due to inaccurately known head and skull shapes and electrode locations in stroke imaging by three-dimensional EIT.

\subsection{Approximation error model for geometric uncertainties in stroke~EIT}

Let 
\begin{equation}\label{addit_accurate} 
  \Vmeas = \Uacc + e,
\end{equation}
denote an accurate enough discretized model for the dependence between the unknowns and the measurements, that is, the associated discretization error is assumed to be negligible. As mentioned above, $\vartheta$ contains the shape parameters $\alpha$ defining the anatomic layers via \eqref{eq:theparam_l} as well as the parameters determining the electrode positions, $\sigbar$ denotes the unknown distributed conductivity in the domain $\domain = \Omega(\vartheta)$, $z$ carries the contact resistances and $e$ models the random measurement errors. In particular, the geometric parameters $\vartheta$ are assumed to accurately define the three-layer anatomic structure for the examined patient as well as the positions of the true physical electrodes of known shape and size.

In practical emergency care, one would usually lack accurate information on the shape of the patient's head $\domain$ and the electrode locations. In such a case, the image reconstruction has to be carried out using an approximate model domain $\modeldomain$, which could be selected, e.g., as an average head from an anatomical atlas with intended electrode locations. In the present setup, this choice corresponds to setting $\alpha = 0$ in the principal component model \eqref{eq:theparam_l} for the three-layer anatomical structure and selecting the parameters defining the electrode positions appropriately; as mentioned above, we denote this choice of geometric parameters by $\vartheta_0$. In other words, the accurate measurement model (\ref{addit_accurate}) is approximated by 
\begin{equation}\label{addit_reduced}
  \Vmeas \approx \mathcal{U}(\sig,z,\vartheta_0) + e = \Umod + e,
\end{equation}
where 
we have used the shorthand notation introduced in \eqref{eq:inv_forv}, $\sig$ is the conductivity in $\modeldomain$ and $z$ still denotes the contact resistances at the $M$ electrodes. The inaccurate forward model $\Umod$ is thus employed as such in the solution of the inverse problem of recovering a conductivity in $\modeldomain$. 

In the approximation error approach, instead of directly utilizing the crude approximation (\ref{addit_reduced}), the {\em accurate measurement model}
(\ref{addit_accurate}) is written in the form
\begin{align}
  \Vmeas &= \, \Umod + \big(\Uacc - \Umod\big) +e \nonumber \\[1mm]
   &= \, \Umod + \varepsilon (\sigbar,\sig,z,\vartheta) + e \nonumber \\[1mm]
&=: \, \Umod + \eta,  \label{addit_accurate2}
\end{align}
where we denote $\eta = \varepsilon + e$ and
$\varepsilon = \varepsilon(\sigbar,\sig,z,\vartheta)$ is the modeling error 
due to the incorrect boundary $\partial \modeldomain$, electrode locations and
skull shape in the forward model.

The objective in the approximation error approach is to derive a model for the likelihood density $\pi (\Vmeas \, \vert \, \sig,z)$ based on the observation model \eqref{addit_accurate2}; here and in what follows, $\pi$ is a generic finite-dimensional probability density whose exact meaning should be clear from the context. Formally, this is achieved via marginalization of the joint likelihood 
$$ 
\pi (\Vmeas \, \vert \, \sig,z ) = \int \int \pi (\Vmeas, \varepsilon, e  \, \vert \, \sig,z) \, {\rm d}e \, {\rm d} \varepsilon
$$ 
with respect to the uninteresting and uncertain parameters in the model \cite{kaipio2013}. In a purely linear and Gaussian case, the marginalization has a closed form solution. In nonlinear and/or non-Gaussian cases, a computationally useful closed form approximation for $\pi (\Vmeas \, \vert \, \sig, z)$ can be obtained by approximating the joint density of the primary unknowns and the nuisance parameters by a Gaussian. 

Let us denote $\omega = (\sig, z)$. Following \cite{kolehmainen2011,kaipio2013} and utilizing a Gaussian surrogate for $\pi(\eta,\omega)$, the approximation for the likelihood density becomes
\begin{align}
\pi(\Vmeas \, \vert \, \omega ) & \varpropto\,\exp \! \left(  
-\frac{1}{2}\big(\Vmeas - \mathcal{U}_0(\omega) - \etamean\big)^{\rm T} 
\etacov^{-1}
  \big(\Vmeas - \mathcal{U}_0(\omega) - \etamean \big)\right),
\label{fem-likelihood} 
\end{align}
where 
{\setlength{\arraycolsep}{0.0em}
\begin{eqnarray}
\etamean&\;=\;& \eta_\ast + \Gamma_{\eta \omega} \Gamma_{\omega}^{-1} 
(\omega - \omega_*), \label{emodel1} \\[1mm]  
\etacov&\;=\;&\Gamma_{\eta} - \Gamma_{\eta \omega}  
\Gamma_{\omega}^{-1} \Gamma_{\eta \omega}^{\rm T} \label{emodel2} .
\end{eqnarray}
}%
Here $\eta_\ast = e_\ast + \epsmean$ is the sum of the means of the measurement and approximation errors, $\Gamma_\eta$ is the covariance of their sum, and $\omega_\ast = (\sigmean, z_\ast)$ and $\Gamma_\omega$ are the mean and covariance of the primary unknown, respectively. Moreover, $\Gamma_{\eta \omega} = \Gamma_{\varepsilon \omega} +\Gamma_{e \omega}$, with $\Gamma_{\varepsilon \omega}$ and $\Gamma_{e \omega}$ being the cross-covariance matrices of the variable pairs $\varepsilon$, $\omega$ and $e$, $\omega$, respectively. In the standard setting, the measurement errors $e$ and the approximation error $\varepsilon$ as well as  $e$ and the primary unknown $\omega$ are assumed to be mutually independent, which means that $\Gamma_{\eta} = \Gamma_{\varepsilon} +\Gamma_{e}$ and $\Gamma_{\eta \omega} = \Gamma_{\varepsilon \omega}$ in \eqref{emodel1}--\eqref{emodel2}.

While it is clear that $\varepsilon$ and $\omega$ are 
not independent, it has regardless turned out in several applications that a feasible model is obtained by setting $\Gamma_{\varepsilon\omega} = 0$. Utilizing this approximation and the earlier assumption that $e$ and $\omega$ as well as $e$ and $\varepsilon$ are mutually independent,
we have
\begin{equation} \label{enhmod}  
\etamean \approx \varepsilon_\ast + e_\ast, \quad 
\etacov \approx \Gamma_{\varepsilon} + \Gamma_{e}
\end{equation}  
in \eqref{emodel1}--\eqref{emodel2}. This approximation, called the 
{\em enhanced error model}
\cite{kaipio04,kaipio07}, is employed in this paper.

\subsection{Learning based estimates for the error statistics}
\label{ssec:estimates}
If the measurement model is linear and the prior model is Gaussian, the approximation error statistics $\epsmean$ and $\epscov$
can be computed analytically,
see \cite{kaipio04}. In nonlinear cases, such as EIT, these
statistics have to be estimated by a simple machine learning procedure.

In our case, the learning algorithm consists of the following four steps:
\begin{enumerate}
\item Generate an ensemble of $N_s$ geometric parameter vectors $\vartheta^{(\ell)}$, each of which characterizes a three-layer head anatomy and an electrode configuration on its exterior surface. The shape parameters are drawn utilizing the principal component based prior model \eqref{eq:theparam_l}--\eqref{eq:Gamma_alpha3}, whereas the spherical coordinates for each electrode center point are independently drawn from a Gaussian distribution whose expected value corresponds to the intended position for that electrode. For each sample configuration, draw random contact resistances $z^{(\ell)}$ and generate a piecewise constant conductivity $\sigbarsamp{\ell}$ by drawing random conductivity levels for the scalp, skull and brain tissues. Solve the {\em accurate} forward model
$\Uaccsamp{\ell}$ for each of the $N_s$ samples $(\sigbarsamp{\ell}, z^{(\ell)}, \vartheta^{(\ell)})$. (The precise probability models for the electrode locations, the contact resistances and the conductivities of the different tissue types are described in Section~\ref{sec:numer}).

\item 
Compute the reference solutions $\Umodsamp{\ell}$ for the $N_s$ pairs $(\sigmean, z^{(\ell)})$, where $\sigmean$ is always the same piecewise constant coefficient corresponding to literature values for the conductivities of the considered three tissue types.

\item Compute the realizations 
\begin{equation}
\label{eq:realizations}
\epssamp{\ell} = \Uaccsamp{\ell} - \Umodsamp{\ell} , \qquad \ell = 1 , \ldots, N_s,
\end{equation}
of the approximation error. Estimate its second order statistics as
\begin{equation}
\label{eq:statistics}
\epsmean = \frac{1}{N_s} \sum_{\ell=1}^{N_s} \epssamp{\ell}, \qquad  
\epscov =  \frac{1}{N_s-1} \sum_{\ell=1}^{N_s} (\epssamp{\ell} - \epsmean) (\epssamp{\ell} - \epsmean)\transp.
\end{equation}
\item Form the enhanced error model via \eqref{fem-likelihood} and \eqref{enhmod} assuming the second order statistics of the measurement noise $e$ are also known. 
\end{enumerate}

\begin{remark} 
The estimation of the approximation error statistics requires altogether $2N_s$ forward solutions in steps 1 and 2 of the above procedure. This can be computationally quite expensive in three-dimensional EIT. However, the statistics only need to be estimated once for a given measurement paradigm, and the associated computations can also be done offline, i.e.~prior to having the measurements in hand (cf.~Section \ref{ssec:simulations}).
In the actual application, the image reconstruction with the approximation error model has a similar computational complexity as the conventional measurement error model \eqref{addit_reduced}.
\end{remark}

\section{Edge-enhancing reconstruction algorithm}
\label{sec:algo}
The potentials measured at the electrodes are modeled as
\begin{equation}
\label{noise_model}
\mathcal{V} = \mathcal{U}_0(\sigmean + \kappa,z) + \eta,
\end{equation}
where $\eta \in \R^{M(M-1)}$ is a realization of a Gaussian random variable with zero mean\footnote{In fact, due to the approximation error approach, the noise does not have a vanishing mean in our numerical experiments; see Section~\ref{sec:approx}. However, this is compensated by abusing the notation and absorbing the mean into the measurement $\mathcal{V}$.} and a known symmetric positive definite covariance matrix $\Gamma \in \R^{M(M-1)\times M(M-1)}$, and the discretized average forward model is as defined in \eqref{eq:inv_forv}. The background conductivity $\sigmean\in \R_+^{N_0}$ corresponds to prior information on (i.e.,~literature values for) the expected conductivities of the considered three tissue types; in our simple three-layer model, $\sigmean$ is assumed to be piecewise constant and compatible with the mean anatomy of the principal component model \eqref{eq:theparam_l}. Our hypothesis is that reconstructing the conductivity perturbation $\kappa\in \R^{N_0}$ in this `average model' produces useful information on the stroke in the true measurement geometry, if the approximation error noise is included in \eqref{noise_model}. See \eqref{fem-likelihood} and \eqref{enhmod} in the preceding section for information on how one arrives at the measurement model \eqref{noise_model} by resorting to enhanced error modeling, but also note that the analysis in this section does not depend on the origin of the additive Gaussian noise term in \eqref{noise_model}.  

We denote a Cholesky decomposition for the inverse of the noise covariance as $\etacov^{-1} = G^{\rm T}G$. The likelihood function,~i.e.~the probability density of the measurement $\mathcal{V}$ given the parameters $\kappa$ and $z$, can be written as (cf.~\eqref{fem-likelihood})
$$
\pi(\mathcal{V} \, | \, \kappa, z) \, \propto \,
\exp \! \Big(-\frac{1}{2} \big| G \big(\mathcal{V} - \mathcal{U}_0(\sigmean + \kappa ,z)\big)\big|^2\Big).
$$
In our hybrid approach, we introduce regularization for the (discretized) conductivity perturbation $\kappa \in \R^{N_0}$ via a prior density 
\begin{equation}
\label{eq:sigma_prior}
\pi(\kappa) \propto \exp \! \big( -\gamma  R(\kappa) \big),
\end{equation}
where $\gamma > 0$ is a free parameter independent of $\kappa$, and $R$ is of the form
\begin{equation}
\label{eq:aRRa}
R(\kappa) = \int_{\Omega_0} \upsilon(x) \, r\big(|\nabla \kappa(x) | \big) \, {\rm d} x.
\end{equation}
Here $r: \R \to \R_+$ corresponds to the prior information that $\kappa$ is approximately piecewise constant, whereas $\upsilon: \Omega_0 \to \R_+$ is the reciprocal of a smooth cut-off function,
\begin{equation}
    \label{eq:upsilon}
\upsilon(x) = \Big( \frac{1}{2}\big( 1 + \tanh( c_\upsilon  ({\rm dist}(x, \partial \Omega_0) - d_\upsilon) ) \big) \Big)^{-1}, \qquad x \in \Omega_0,
\end{equation}
where $c_\upsilon > 0$ and $d_\upsilon > 0$ are to be chosen so that the value of $\upsilon$ is large in the skin layer and decreases to almost one within the skull layer. The inclusion of $\upsilon$ in \eqref{eq:aRRa} reflects the prior information that the stroke is expected to only affect the conductivity in the brain. Moreover, it is reasonable to assume the skull has more local variations in its conductivity than the skin. (Recall also that variations in the base conductivity levels for the skin and skull layers are already accounted for in the approximation error model.) We exclusively employ a smoothened TV prior that corresponds to the choice~\cite{Rudin92}
\begin{equation}
\label{eq:arra}
r(t) = \sqrt{T^2 + t^2} \approx | t |
\end{equation}
but one could as well consider,~e.g.,~a PM or smoothened TV$^q$ prior. Here, $T>0$ is a small parameter that ensures the differentiability of $r$. The contact resistances are given an {\em uninformative} prior, that is, all realizations of $z$ in $\R^M$ are considered to be equally probable {\em a priori}.\footnote{However, in the actual reconstruction algorithm the components of $\sigma$ and $z$ are clamped onto physically reasonable subintervals of $\R_+$ to ensure stable solvability of forward problems.}

According to the Bayes' formula, the posterior density for $(\kappa, z)$ is
\begin{align*}
\pi(\kappa, z \, | \, \mathcal{V}) \, &\propto \, \pi(\mathcal{V} \, | \, \kappa, z) \, \pi(\kappa)  \\ 
&\propto \, \exp \Big(-\frac{1}{2} \big| G \big(\mathcal{V} - \mathcal{U}_0(\sigmean + \kappa ,z)\big)\big|^2 - \gamma R(\kappa) \Big).
\end{align*}
Hence, finding a MAP estimate for the pair $(\kappa, z)$ corresponds to minimizing the Tikhonov functional 
\begin{equation}
\label{eq:Tikhonov}
\Phi(\kappa,z) :=  \frac{1}{2} \big| G \big(\mathcal{V} - \mathcal{U}_0(\sigmean + \kappa ,z)\big)\big|^2 + \gamma R(\kappa).
\end{equation}
The algorithm introduced in \cite{Harhanen15} is motivated by this task, although it slightly deviates from the strict Bayesian interpretation and utilizes instead efficient Krylov subspace iterations and the Morozov discrepancy principle.

Before writing down the actual algorithm, we still need to introduce a couple of auxiliary concepts; for more information consult~\cite{Harhanen15}. A preliminary $\kappa$-dependent regularization matrix $\tilde{H}(\kappa) \in \R^{N_0 \times N_0}$, encoding information on the variations in $\kappa$, is defined via
\begin{equation}
\label{eq:H}
\tilde{H}_{k,l}(\kappa) := \int_{\Omega_0} \upsilon(x) \frac{r'(|\nabla \kappa(x)|)}{|\nabla \kappa(x)|} \nabla \varphi_k^{(0)}(x) \cdot \nabla \varphi_l^{(0)}(x) \, {\rm d} x, \qquad k,l=1,\dots, N_0.
\end{equation}
Observe that $\tilde{H}(\kappa)$ corresponds to a FE discretization in the basis $\{ \varphi_k^{(0)} \}_{k=1}^{N_0} \subset H^1(\Omega_0)$ for the partial differential operator
\begin{equation}
\label{eq:diffop}
- \nabla \cdot \big( \beta (|\nabla \kappa|)  \nabla (\, \cdot \,) \big), \qquad \text{with} \ \ \beta(w) := \upsilon \frac{r'(w)}{w}
\end{equation}
and a homogeneous Neumann condition on $\partial \Omega_0$. It is straightforward to check that
$$
c \leq \beta(w) \leq C, \qquad w \in L^\infty(\Omega_0),
$$
 almost everywhere in $\Omega_0$ for $c, C > 0$, of which $C$ only depends on $T$ and $\upsilon$ and $c$ only on $T$ an upper bound for $\| w \|_{L^\infty(\Omega_0)}$. It follows that the operator \eqref{eq:diffop} is uniformly elliptic in $\kappa$ over any bounded subset of $\R^{N_0}$. In particular, $\tilde{H}(\kappa)$ is positive semi-definite with $[1, \dots, 1]^{\rm T} \in \R^{N_0}$ spanning its kernel.

To ensure the invertibility of the actual regularization matrix employed in the reconstruction algorithm, we define it through
\begin{equation}
\label{eq:H2}
H(\kappa) = \tilde{H}(\kappa) + \lambda(\kappa) \mathbb{I},
\end{equation}
where $\lambda(\kappa)$ is the smallest nonzero eigenvalue of $\tilde{H}(\kappa)$ and $\mathbb{I} \in \R^{N_0 \times N_0}$ denotes the identity matrix. We write $H(\kappa) = L(\kappa)^T L(\kappa)$ for a Cholesky decomposition of $H(\kappa)$, although such a factorization need not be explicitly computed at any stage of the following reconstruction algorithm.

\begin{algorithm}
\label{alg:kokohoska}
\vspace{2mm}
Pick $T > 0$ for $r: \Omega_0 \to \R_+$ in \eqref{eq:arra}, $c_\upsilon$ and $d_\upsilon$ for \eqref{eq:upsilon}, the number of lagged diffusivity steps $N_{\rm LD} \in \N$ for the interior loop, and the feasibility intervals $[\kappa_{\rm min}, \kappa_{\rm max}]$, $[z_{\rm min}, z_{\rm max}] \subset \R_+$ for the components of $\kappa$ and $z$, respectively. Set $\kappa^{(0)} = 0$ and choose $z^{(0)} \in \R^M_+$ to be the minimizer of
\begin{equation}
\label{residual}
E(\kappa^{(0)}, z) := \big | G \big(\mathcal{V} - \mathcal{U}_0(\sigmean + \kappa^{(0)} ,z)\big) \big|^2
\end{equation}
with respect to $z$ in the span of $[1, \dots, 1]^{\rm T} \in \R^M$. 
Set $j=0$, $\epsilon = \sqrt{M(M-1)}$ and $\mathcal{U}_0 = \mathcal{U}_0(\sigmean+\kappa^{(0)},z^{(0)}) = \mathcal{U}_0(\sigmean,z^{(0)})$.  
\vspace{2mm}
\begin{enumerate}
\item Evaluate the Jacobian matrices of the map $(\sigma, z) \mapsto \mathcal{U}_0(\sigma, z) $ with respect to $\sigma$ and $z$ at $(\sigmean + \kappa^{(j)}, z^{(j)})$. Denote them by $J_1 \in \R^{M(M-1) \times N_0}$ and $J_2 \in \R^{M(M-1) \times M}$, respectively.
\vspace{2mm}
\item Set $y = G(\mathcal{V} - \mathcal{U}_0 + J_1\kappa^{(j)} + J_2 z^{(j)})$, $B_1 = G J_1$ and $B_2 = G J_2$.  
\vspace{2mm}
\item Let $Q$ be the orthogonal projection onto $\mathcal{R}(B_2)^{\perp}$, where $\mathcal{R}(B_2)$ denotes the range of the matrix $B_2$. Define $A = Q B_1$ and $b = Qy$.
\vspace{2mm}
\item Set $l=0$, $\kappa^{(j,0)} = \kappa^{(j)}$ and $z^{(j,0)} = z^{(j)}$.
For $l=0, \dots, N_{\rm LD}-1$: 
\begin{enumerate}
\vspace{2mm}
\item Build the preconditioner $H = H(\kappa^{(j,l)})$ according to \eqref{eq:H} and \eqref{eq:H2}. 
\vspace{2mm}
\item Apply the {\em LSQR} algorithm from \cite{Arridge13} to 
$$
(L^{-1})^{\rm T} A^{\rm T} A L^{-1}  \tilde{\kappa} \, = \,  (L^{-1})^{\rm T} A^{\rm T} b, \qquad \kappa = L^{-1} \tilde{\kappa},
$$
starting from  $\kappa = 0$, which results in a solution sequence $\{ \kappa^{(j,l)}_k\}_{k\geq 0}$. Stop the iteration when $\| A \kappa^{(j,l)}_k - b \| \leq \epsilon$. Denote the resulting solution as $\kappa^{(j,l+1)}$.
(Recall that $L$ corresponds to a symmetric decomposition of the form $H = L^{\rm T} L$, but it need not be formed explicitly \cite{Arridge13}.)
\vspace{2mm}
\item Set
$z^{(j,l+1)} = (B_2^{\rm T} B_2 )^{-1} B_2^{\rm T} (y  - B_1 \kappa^{(j,l+1)})$.
\vspace{2mm}
\item Clamp the components of $\kappa^{(j,l+1)}$ and $z^{(j,l+1)}$ onto intervals $[\kappa_{\rm min}, \kappa_{\rm max}]$ and $[z_{\rm min}, z_{\rm max}]$, respectively.
\end{enumerate}
\vspace{2mm}
\item Set $\kappa^{(j+1)} = \kappa^{(j,N_{\rm LD})}$ and $z^{(j+1)} = z^{(j,N_{\rm LD})}$.
\item Compute $\mathcal{U}_0 = \mathcal{U}_0(\sigmean + \kappa^{(j+1)}, z^{(j+1)})$. If $E(\kappa^{(j+1)}, z^{(j+1)}) \leq \epsilon$, dub $(\sigmean + \kappa^{(j+1)}, z^{(j+1)})$ the reconstruction and exit. Otherwise, set $j \leftarrow j+1$ and return to step~1.
\end{enumerate}
\vspace{2mm}
\end{algorithm}

Let us briefly explain the different phases in Algorithm~\ref{alg:kokohoska}; see once again \cite{Harhanen15} for more details. After estimating the homogeneous contact resistances that most accurately explain the data under the assumed noise model and expected anatomy, step~(1) corresponds to the linearization of the forward map around the current estimate for the conductivity and the contact resistances; see,~e.g.,~\cite{Kaipio00,Lechleiter06}. In step~(2), this linearization is plugged in the data fit term of~\eqref{eq:Tikhonov} so that the unknowns are chosen to be the conductivity perturbation $\kappa$ and the contact resistances~$z$, not the additive changes in them, and thus the data vector becomes $y$. Defining $A = Q B_1$ and $b = Q y$ in step~(3) sets, loosely speaking, the aim of only explaining the projection of $y$ onto the orthogonal complement of $\mathcal{R}(B_2)$ by the conductivity perturbation. 

The contribution of the penalty term $R(\kappa)$ to the necessary condition for a minimizer of~\eqref{eq:Tikhonov} can be written as $(\nabla R)(\kappa) = H(\kappa) \kappa$, not accounting for the small extra diagonal term in \eqref{eq:H2}. The loop in step~(4) originates from approximating this term in the spirit of the lagged diffusivity iteration as $H = H(\kappa^{(j,l)}) \kappa$ \cite{Vogel96},~i.e.,~to evaluating the regularization matrix at a previous estimate for $\kappa$. This leads in principle to simple Tikhonov regularization for computing $\kappa^{(j,l+1)}$,
\begin{equation}
\label{eq:tikhonov}
(A^{\rm T}A + \gamma H)\kappa = A^{\rm T} b,
\end{equation}
where $H$ carries information on the variations observed in the previous estimate for~$\kappa$. 

Instead of solving \eqref{eq:tikhonov} as such, the regularization induced by $\gamma>0$ is replaced in step (4b) by an LSQR iteration for a priorconditioned (cf.~\cite{Calvetti07,Calvetti12,Calvetti05,Calvetti08}) version of $A^{\rm T}\!A \kappa = A^{\rm T} b$ accompanied with an early stopping rule motivated by the Morozov discrepancy principle and the assumed noise level. To put it very short, the symmetric preconditioning with $H$ ensures that the conductivity perturbation produced by LSQR always lies in the range of $H^{-1}$, which consists of elements that are in accordance with the prior information from the previous round encoded in $H$. Finally, step~(4c) corresponds to choosing $z^{(j,l+1)}$ to be the solution of $B_2z = y - B_1 \kappa^{(j,l+1)}$ in the sense of least squares; in the linearized framework, one tries to explain by $z^{(j,l+1)}$ everything in the data that was not already explained by the choice of~$\kappa^{(j,l+1)}$.

The whole algorithm is stopped in step (6) if the original nonlinear equation satisfies the Morozov discrepancy principle. For an explanation why $\epsilon = \sqrt{M(M-1)}$ is a reasonable choice for the discrepancy level in both the inner and the outer loop of Algorithm~\ref{alg:kokohoska}, see ones again \cite{Harhanen15}.

\begin{remark}
Algorithm~\ref{alg:kokohoska} differs from its original version presented in \cite{Harhanen15} in a few ways. First of all, in~\cite{Harhanen15} the aim is to reconstruct the conductivity itself, not its deviation from some known conductivity level that in our case corresponds to the the prior mean conductivity of the average head model. In particular, \cite{Harhanen15} assumes that the actual conductivity inside the examined body is {\em a priori} characterized by an edge-promoting prior, whereas here we assume such prior information for the conductivity perturbation $\kappa$. Secondly, \cite{Harhanen15} ensures the invertibility of the regularization matrix by introducing nonhomogeneous Dirichlet conditions for the elliptic differential operator \eqref{eq:diffop} at the electrodes, that is, the conductivity underneath the electrodes is fixed to a certain homogeneous estimate throughout the reconstruction process. In this work, we instead shift the spectrum of the regularization matrix slightly to the right from the origin, as indicated by \eqref{eq:H2}. Thirdly, \cite{Harhanen15} employs exclusively $N_{\rm LD} = 1$, that is, there in only one lagged diffusivity step per a linearization of the forward model in the algorithm of \cite{Harhanen15}. In our numerical experiments, we instead choose $N_{\rm LD} = 5$: due to the relatively high amount of approximation error noise in the data, the exterior loop of Algorithm~\ref{alg:kokohoska} often reaches the Morozov level after only a couple of linearizations of the forward model, which does not give enough time for the lagged diffusivity steps to sufficiently enhance inclusion boundaries unless more than one step per a linearization of the forward model is allowed. Finally, \cite{Harhanen15} does not consider spatially varying priors, as is \eqref{eq:sigma_prior} due to the weight function $\upsilon$ in~\eqref{eq:aRRa}.
\end{remark}

\section{Numerical examples}
\label{sec:numer}
This section presents our numerical examples.
We start by briefly explaining how realistic measurement data are simulated and how the approximation error statistics are computed. Thereafter, some numerical results are presented with emphasis on comparing reconstructions obtained with and without adopting the approximation error approach for different levels of geometric mismodeling in the measurement setup.    

\subsection{Simulation of approximation error statistics and measurement data}
\label{ssec:simulations}
\begin{figure}
\center{
 {\includegraphics[width=6cm]{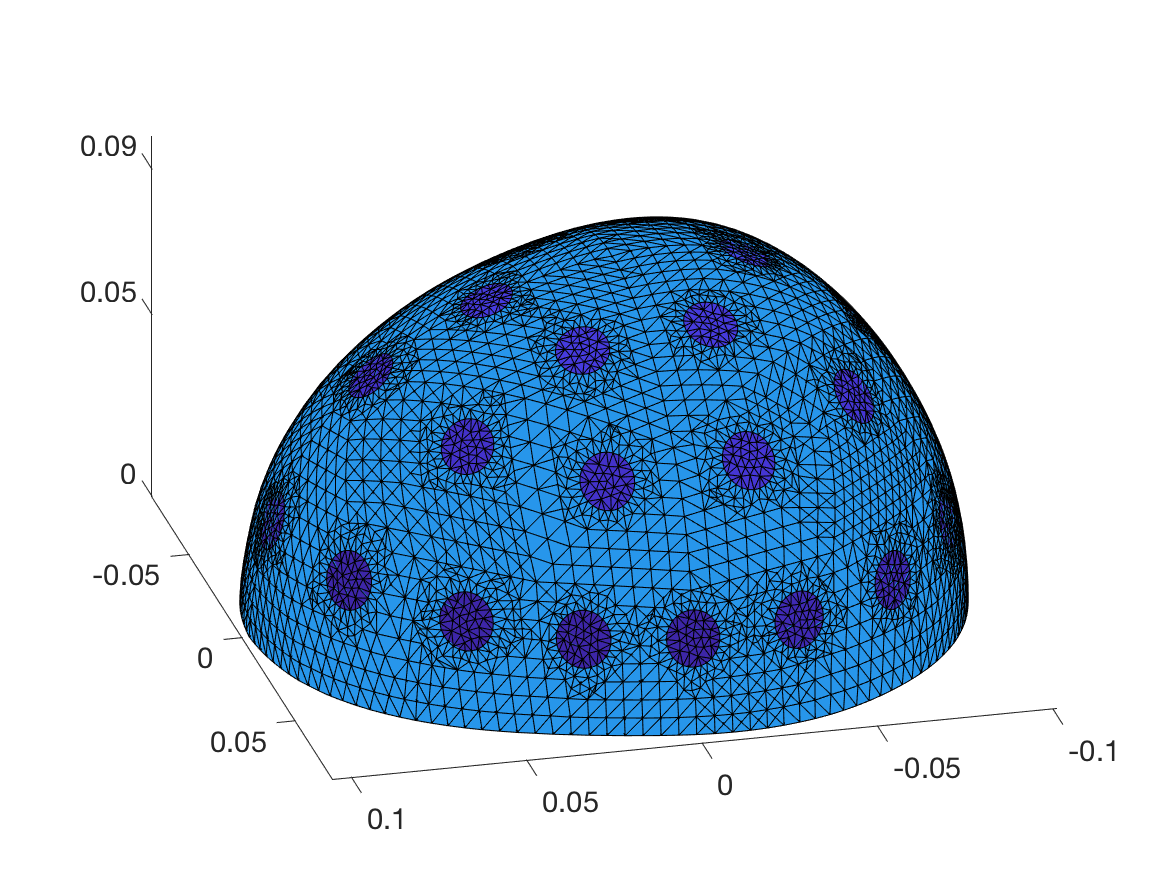}}
  \quad
 {\includegraphics[width=6cm]{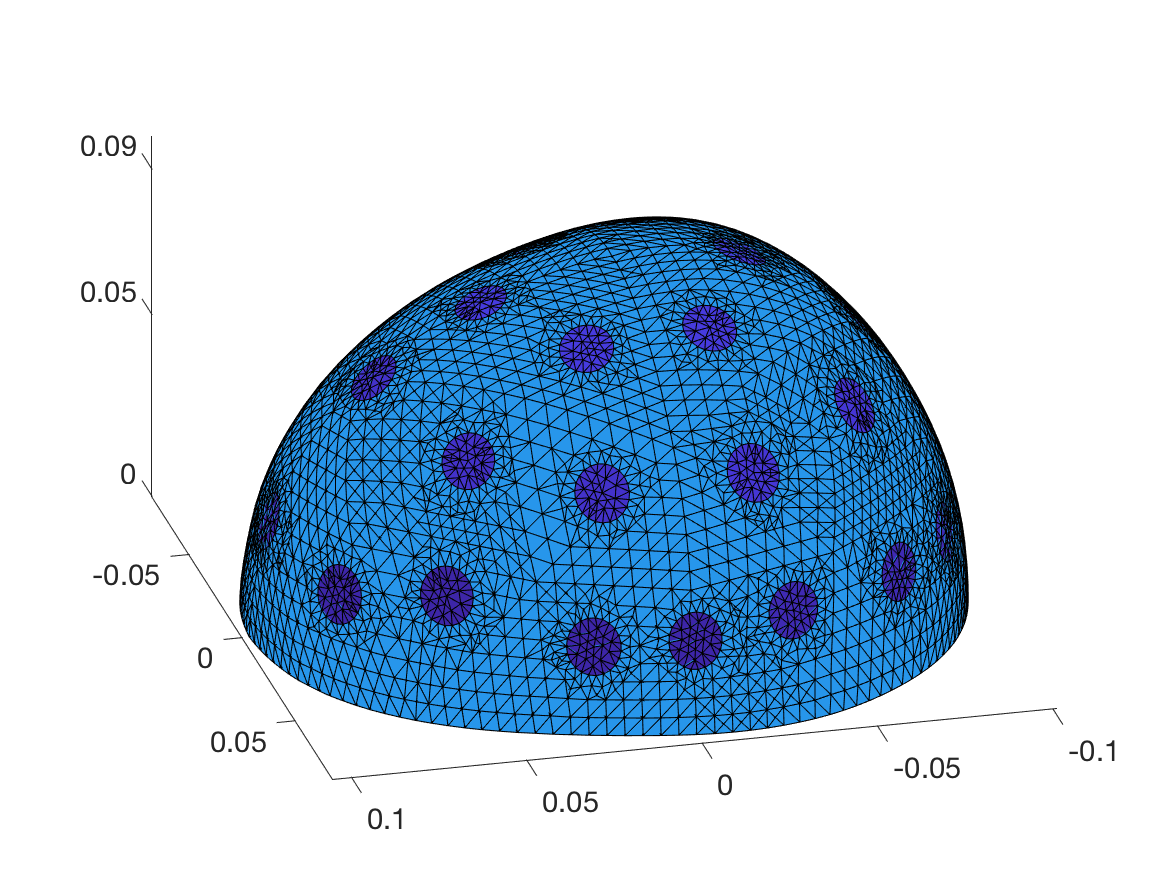}}
 }
\caption{Left: a computational head model with three belts of electrodes at their intended spherical angles $\bar{\theta}$ and $\bar{\phi}$. Right: the same head model with perturbed electrode positions, with $\varsigma_\theta = \varsigma_\phi = 0.015$ in \eqref{eq:Gamma_thph}. 
}
\label{fig:elec}
\end{figure}

All measurements are performed with $M=32$ electrodes that are organized into three horizontal belts around the head as illustrated in Figure~\ref{fig:elec}. The employed current patterns are of the form $I(j)=\text{e}_p-\text{e}_j$, $j=1,\dots,p-1, p+1,\dots, M$, where $p\in  \{ 1, \dots ,M\}$  is the label of the so-called current-feeding electrode, which in our numerical experiments is $E_{27}$,~i.e.~the frontal one on the top electrode belt. The common radius of the electrodes $R=0.75\,{\rm cm}$ is assumed to be known.
When estimating the approximation error statistics in the manner explained in Section~\ref{ssec:estimates}, the electrode potential measurements corresponding to the accurate model $\Uacc$ and the surrogate model $\mathcal{U}_0(\sigmean,z)$ are simulated via solving \eqref{eq:cemeqs} by a FEM with piecewise linear basis functions on a dense mesh (cf.~\cite{Vauhkonen97, Vauhkonen99}). Subsequently, the actual statistics are approximated via the steps \eqref{eq:realizations}--\eqref{eq:statistics}. 

As explained in Section~\ref{sec:approx}, we adopt for the surrogate forward model the average three-layer head geometry corresponding to vanishing shape parameters $\alpha_i=0$, $i=1, \dots, \tilde{n} := 10$, and assume the expected values of the polar and azimuthal angles for the electrode centers. These angles, denoted by $\bar{\theta} \in (0, \pi/2)^{M}$ and $\bar{\phi} \in [0, 2\pi)^{M}$, correspond to the intended positions of the electrodes,~i.e.~the angles at which the electrodes should be if they had been positioned properly; see the left-hand image in Figure~\ref{fig:elec}. 

For each realization of the accurate forward model, the shape parameter vector $\alpha \in \R^{\tilde{n}}$ is drawn from the normal distribution $\mathcal{N}(0, \Gamma_\alpha)$, where $\Gamma_\alpha \in \R^{\tilde{n} \times \tilde{n}}$ is the diagonal covariance matrix induced by the employed head library and defined elementwise in \eqref{eq:Gamma_alpha3}. The associated central angles for the electrodes, $\theta$ and $\phi$, are drawn from the distributions $\mathcal{N}(\bar{\theta}, \Gamma_\theta)$ and $\mathcal{N}(\bar{\phi}, \Gamma_\phi)$, with
\begin{equation}
\label{eq:Gamma_thph}
\Gamma_\theta = \varsigma_\theta^2 \mathbb{I} \qquad {\rm and} \qquad
\Gamma_\phi = \varsigma_\phi^2 \mathbb{I}. 
\end{equation}
Here, $\mathbb{I} \in \R^{M \times M}$ is the identity matrix and $\varsigma_\theta,\varsigma_\phi>0$ determine the standard deviations in the two angular directions. Notice that $\varsigma_\theta$ and $\varsigma_\phi$ must be chosen so small that the electrodes are not at a risk to overlap or end up outside the crown of the computational head; in our numerical tests, we use the standard deviations $\varsigma_\theta=\varsigma_\phi = 0.015$ that correspond to an average electrode displacement of about $1.5$--$2$\,mm depending on the height of the electrode belt; see the right-hand image of Figure~\ref{fig:elec}.

As explained in Section~\ref{ssec:estimates}, the contact resistances are common for the accurate and surrogate forward models when the approximation error statistics are simulated. These values, together with the conductivity levels for the three tissue types in the accurate model simulations, are drawn from the associated normal distributions listed in Table~\ref{tab:values}. These distributions are in line with the conductivity levels reported in the medical literature for the considered tissue types~\cite{Gabriel96,Lai05,Latikka01,mccann2019variation,Oostendorp00}. On the other hand, the conductivity levels related to the surrogate model are fixed and correspond to the mean values of the distributions listed in Table~\ref{tab:values}.

\begin{table}[b!]
\centering
\caption{The probability distributions used for different parameters in the simulation of the approximation error statistics.}
\begin{tabular}{|l|l|}
\hline
conductivity of the scalp $\sigma_{\rm scalp}$ & $\mathcal{N}(0.2, 0.02^2) $ (S/m) \\ \hline
conductivity of the skull $\sigma_{\rm skull} $ &  $\mathcal{N}(0.06, 0.006^2) $  (S/m) \\ \hline
conductivity of the brain $\sigma_{\rm brain}$ &  $\mathcal{N}(0.2, 0.02^2) $ (S/m)  \\ \hline
contact resistances $z$ &  $\mathcal{N}(0.01, 0.0025^2 \mathbb{I}) $ ($\Omega \, {\rm m}^2$)  \\ \hline
electrode angles $\theta \in (0,\pi/2)^M$, $\phi \in \R^M$ & $\mathcal{N}(\bar{\theta},0.015^2 \mathbb{I})$, $\mathcal{N}(\bar{\phi},0.015^2 \mathbb{I})$ (rad) \\ \hline
shape parameters $\alpha$ & $\mathcal{N}(0, \Gamma_\alpha)$ \\ \hline
\end{tabular}
\label{tab:values}
\end{table}

For the estimation of the approximation error statistics according to \eqref{eq:realizations}--\eqref{eq:statistics}, $N_s = 1000$ samples of training data were simulated in parallel on Triton~\cite{Triton}, the Aalto University School of Science high-performance computing cluster, on $100$ different nodes with $10$ processors each. Each forward computation, corresponding to $M-1$ linearly independent electrode current patterns, required about $30$ seconds and $3$\,MB of memory on FE meshes that contained around $25\,000$ nodes and $100\,000$ tetrahedra. The overall computation time for generating the training data did not exceed thirty minutes.

The noisy measurements $\mathcal{V}$ for a particular numerical experiment are simulated as indicated by \eqref{addit_accurate}, with the forward solution $\Uacc$ corresponding to an independently drawn set of shape parameters, electrode positions, contact resistances, and background conductivity levels for the three tissue types (cf.~Table~\ref{tab:values}). In addition, there is typically an embedded inhomogeneity in the brain tissue, representing the to-be-identified stroke and taken into account in the parametrized conductivity~$\sigma$ entering the forward model. The measurement noise $e \in \R^{M(M-1)}$ in \eqref{addit_accurate} is a realization of a zero-mean Gaussian with the diagonal covariance matrix 
\begin{equation}
\Gamma_\eta = \varsigma_\eta^2\big( \max_j(\mathcal{U}_j)- \min_j(\mathcal{U}_j) \big)^2 \mathbb{I} \in \R^{M(M-1) \times M(M-1)},
\end{equation}
where the free parameter $\varsigma_\eta>0$ can be tuned to set the relative noise level. In our numerical experiments we choose it systematically as $\varsigma_\eta = 10^{-3}$. Such a noise model has been used with real-world data,~e.g.,~in \cite{Darde13b}.

\subsection{Reconstructions}
We apply Algorithm~\ref{alg:kokohoska} to three test cases. The corresponding target heads are formed according to the following specifications.
\begin{itemize}
    \item[Case~1:]
    The shape of the patient head is inaccurately known (moderate error). The parameters defining the target head  and the electrode configuration are drawn from the distributions listed in Table~\ref{tab:values} after scaling all standard deviations by 0.5.
    \item[Case~2:] The shape of the patient head is inaccurately known (larger error). The parameters defining the target head and the electrode configuration are drawn from the distributions listed in Table~\ref{tab:values}.
    \item[Case~3:] The shape of the patient head is exactly known. The target head and the electrode configuration are defined by the mean parameter values listed in Table~\ref{tab:values}.
\end{itemize}
Take note that the specifications of Cases~1--3 also define the contact resistances and the background conductivity levels for the three tissue types via the associated random draws, meaning that the error increases in the same proportion.

Case~1 corresponds to a target head that is relatively close to the average head used as the surrogate in our reconstruction algorithm, whereas Case~2 corresponds to the expected variation in the uncertain parameters. The third case is a `sanity check' involving the average head, the electrodes at their intended positions, the mean values for the background conductivity levels and the expected contact resistances.

In each geometric case, we test Algorithm~\ref{alg:kokohoska} with and without approximation error modeling for three types of virtual patients:
\begin{itemize}
    \item[(i)] A healthy patient with no stroke, that is, there is no inhomogeneity inside the examined head. 
    \item[(ii)]  A hemorrhagic stroke causing an inclusion with conductivity $2$\,S/m.
    \item[(iii)] An ischemic stroke causing an inclusion with conductivity $0.02$\,S/m.
\end{itemize}
To allow straightforward comparison between different geometric setups and medical conditions, the stroke is located at the same position in all tests: it is characterized by a ball of radius $2.25$\,cm around $(2,3,3)$\,cm, with the origin at the center of the bottom face of the examined head.

Approximation error modeling can either be taken into account or ignored when Algorithm~\ref{alg:kokohoska} is run for a given set of measurements. In the former case, the additive noise $\eta$ in \eqref{noise_model} is the sum of the (mutually independent) measurement noise $e$ and approximation error noise $\varepsilon$ (minus their mean). In the latter case, the approximation error noise is ignored and $\eta = e$ is used. Take note that the covariance matrix of $\eta$ affects Algorithm~\ref{alg:kokohoska} via its Cholesky factor $G$; in particular, the covariance matrix is larger in the sense of positive definiteness when approximation error noise is included in the analysis.

Our choices for the free parameters in Algorithm~\ref{alg:kokohoska} are: $T = 10^{-6}$, $c_\upsilon = 3\, {\rm cm}^{-1}$, $d_\upsilon = 1$\,cm, $N_{\rm LD} = 5$, $z_{\rm min} =10^{-6}\, \Omega\,{\rm m}^2$, $z_{\rm max} =10\, \Omega\,{\rm m}^2$. Furthermore, the bounds $\kappa_{\rm min}$, $\kappa_{\rm max}$ for the conductivity perturbation are chosen so that the reconstructed conductivity always lies in the interval $[10^{-5}, 10^2]$\,S/m. When approximation error modeling is employed, Algorithm~\ref{alg:kokohoska} is implemented in its standard form,~i.e.~it is run until the stopping criterion of its exterior loop is satisfied.
However, this stopping criterion is not typically reached when approximation error noise is not taken into account, due to the Morozov criterion's assumptions not being fulfilled, and thus Algorithm~\ref{alg:kokohoska} is terminated when the residual monitored in the exterior loop of Algorithm~\ref{alg:kokohoska}, i.e.~$E(\, \cdot \, , \, \cdot \,)$, increases for the first time. 

All reconstructions were formed on finite element meshes with approximately $390\,000$ tetrahedrons and $70\,000$ nodes, with appropriate refinements close to the electrodes. The presented results correspond to a MATLAB-based implementation on a laptop with 8 GB RAM and an Intel CPU having clock speed 2.3\,GHz. With approximation error modeling, the running time of Algorithm~\ref{alg:kokohoska} was about 7--15 minutes depending on the number of required exterior loops for reaching the noise level. Without approximation error modeling, the running times were typically somewhat longer, as it took at least a few exterior loops of Algorithm~\ref{alg:kokohoska} before the residual $E(\, \cdot \, , \, \cdot \,)$ started to increase.

\subsubsection{Case~1:~head shape and other parameters close to the average}

\begin{figure}[t!]
\center{
  {\includegraphics[width=4.1cm]{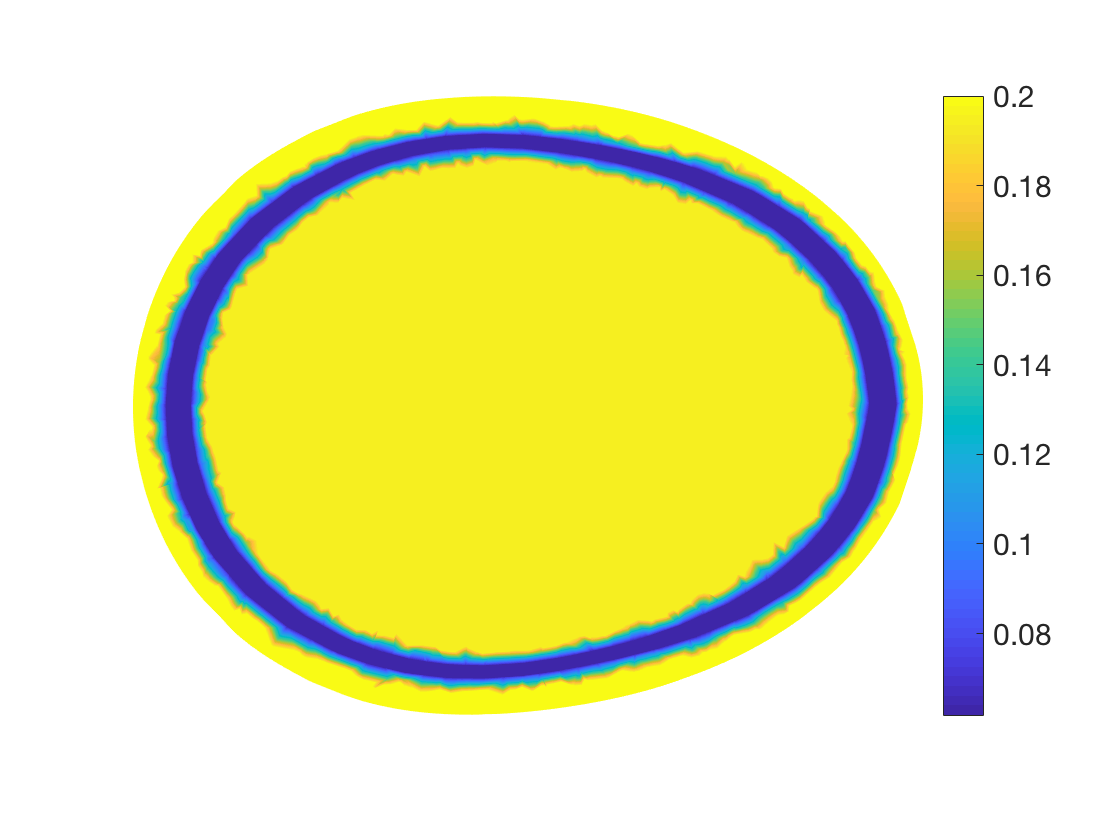}}
  {\includegraphics[width=4.1cm]{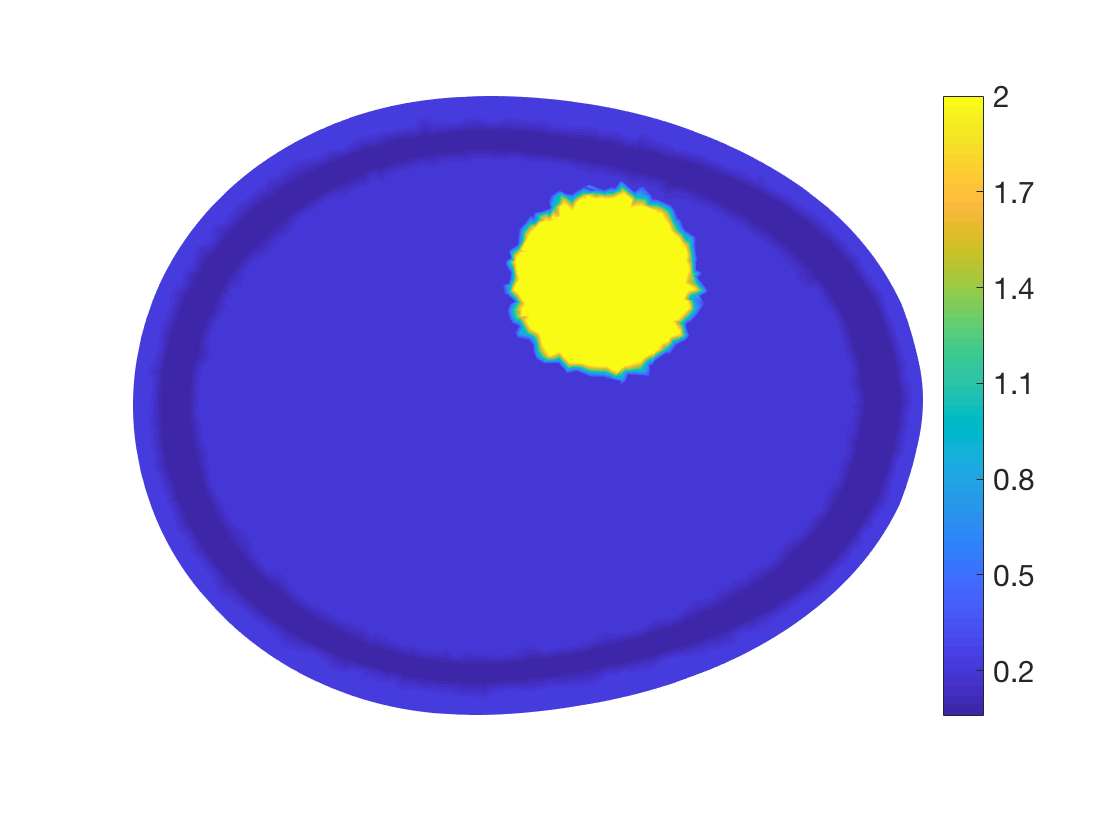}}
  {\includegraphics[width=4.1cm]{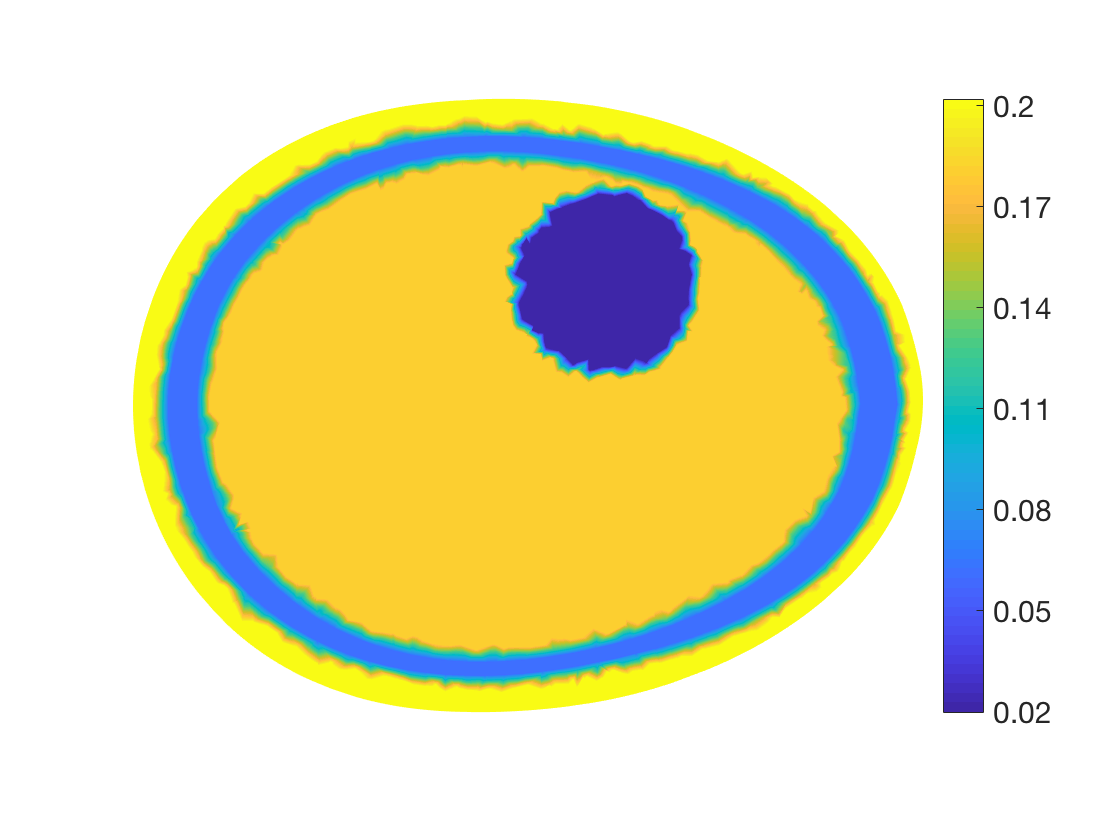}}
  }
  \center{
  {\includegraphics[width=4.1cm]{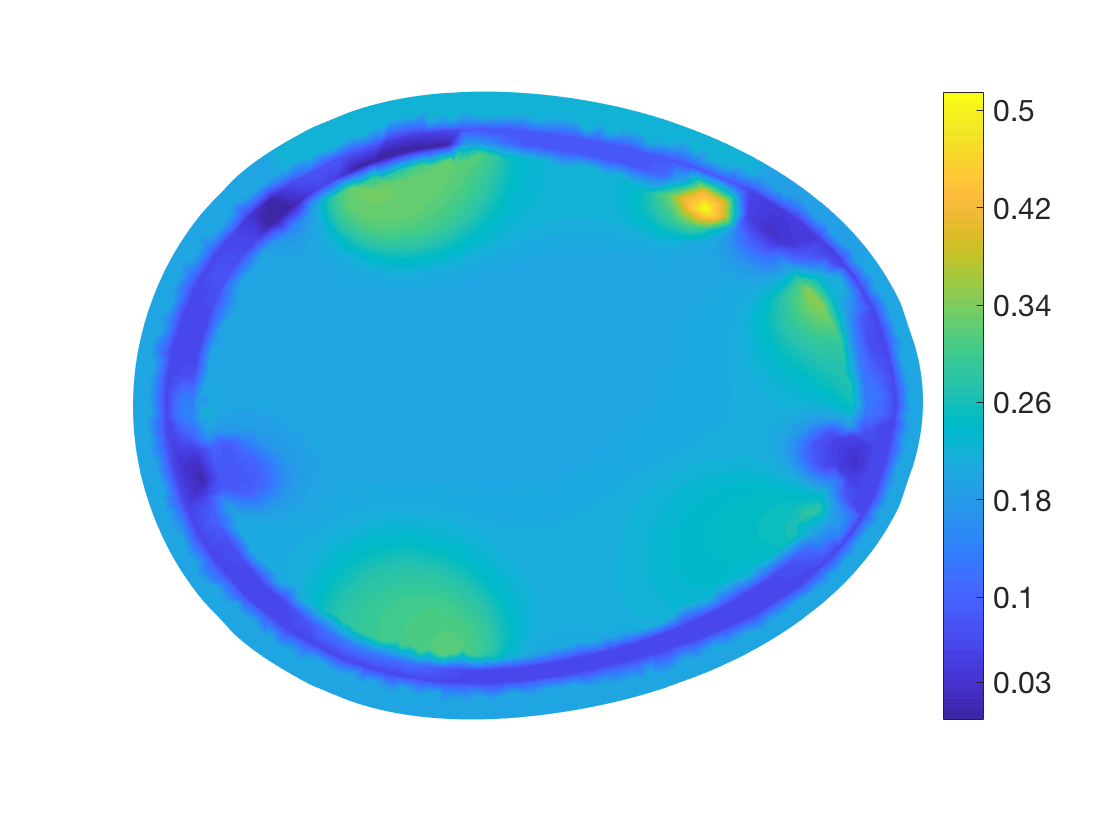}}
  {\includegraphics[width=4.1cm]{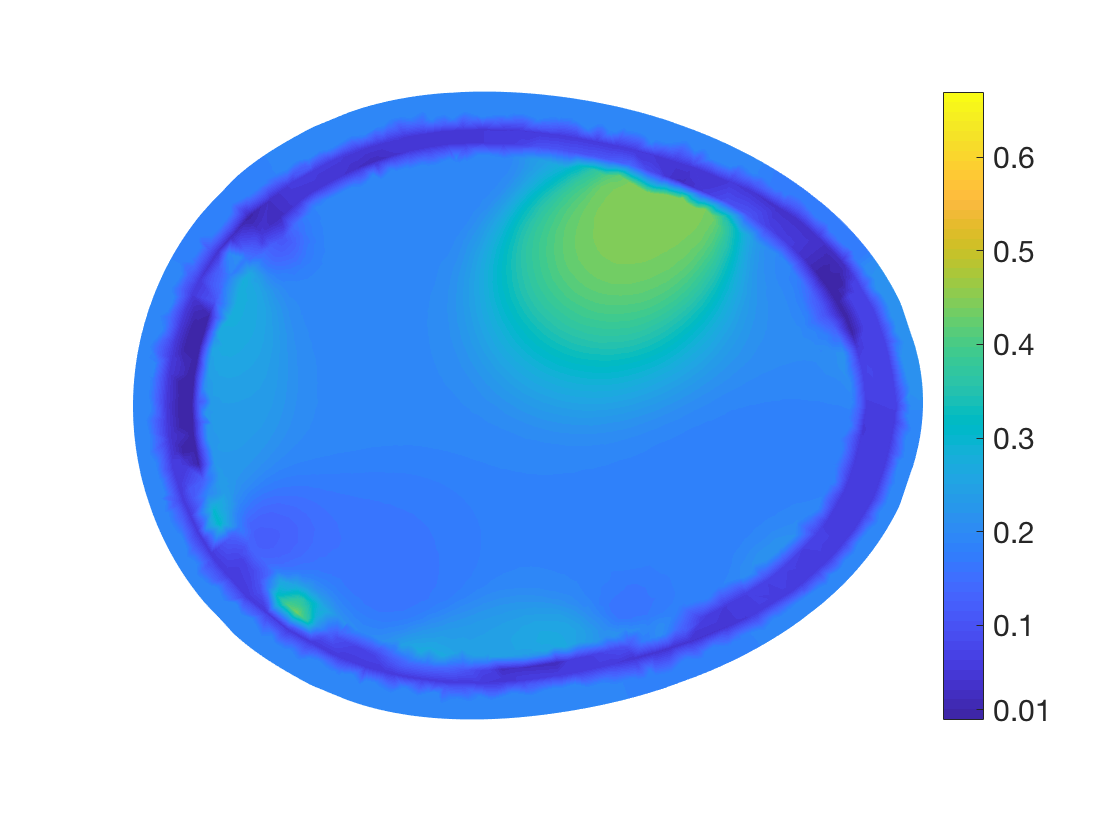}} 
  {\includegraphics[width=4.1cm]{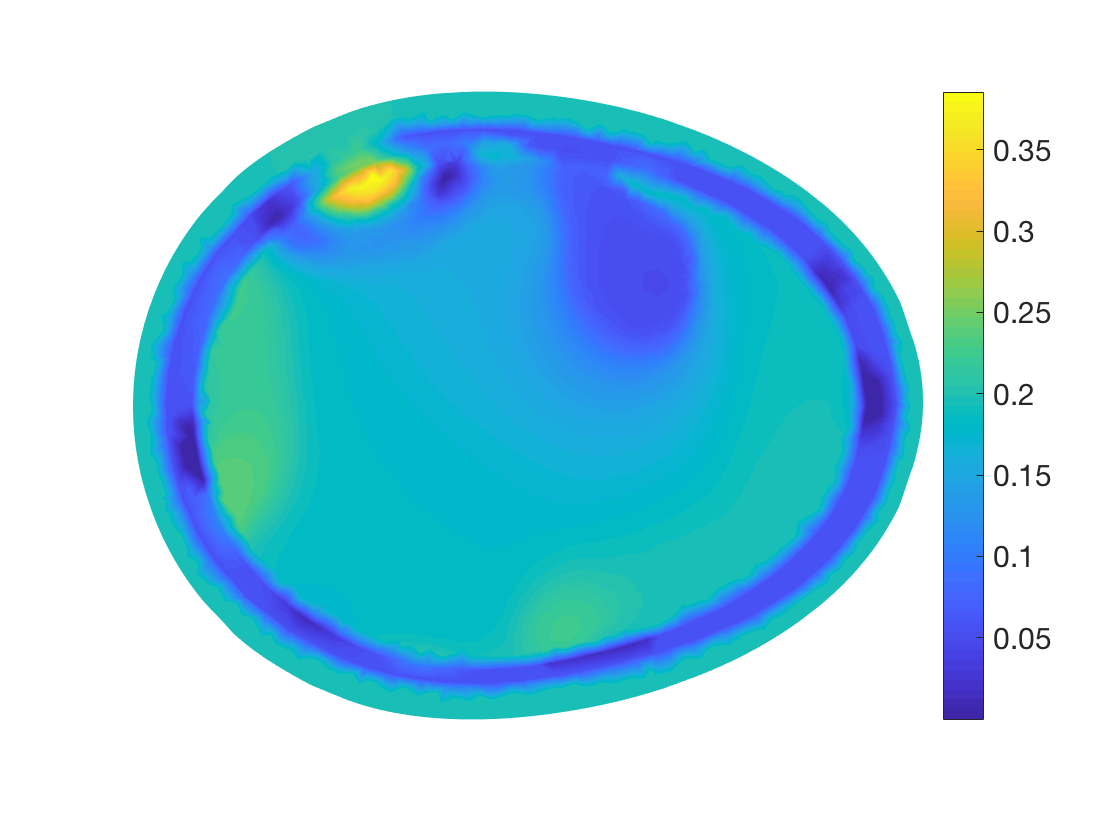}}
  }
  \center{
  {\includegraphics[width=4.1cm]{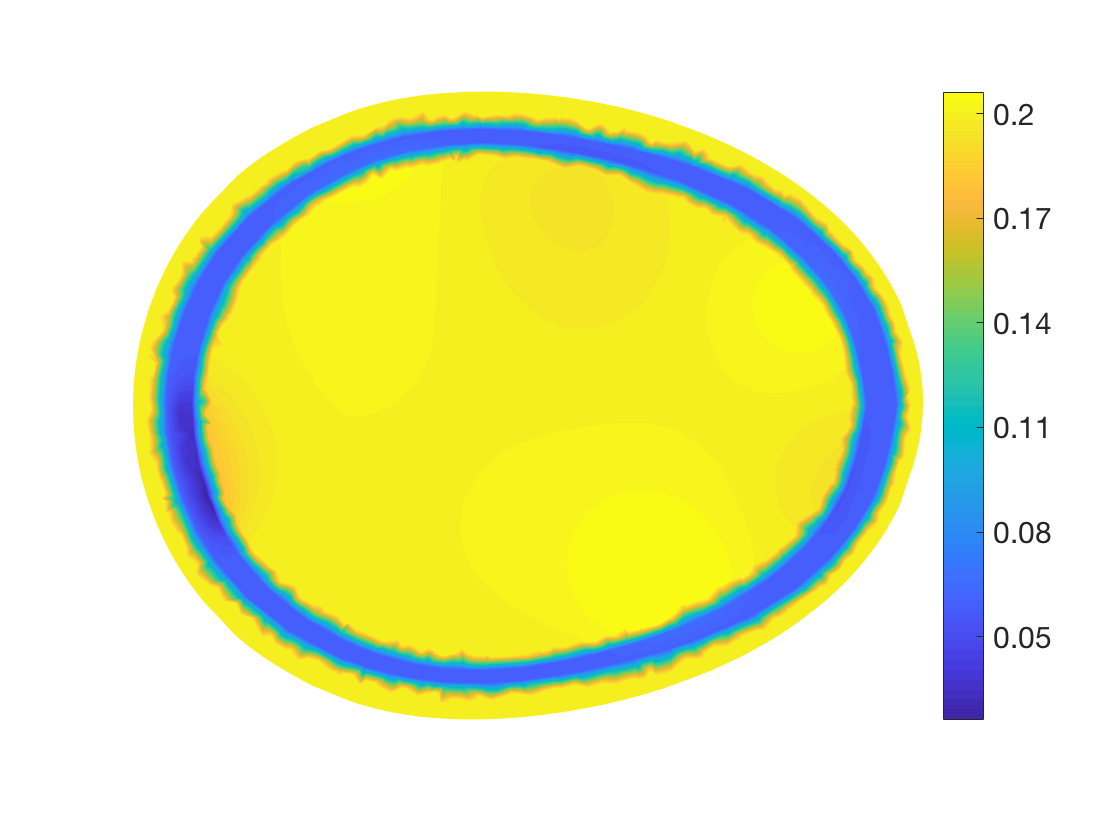}}
  {\includegraphics[width=4.1cm]{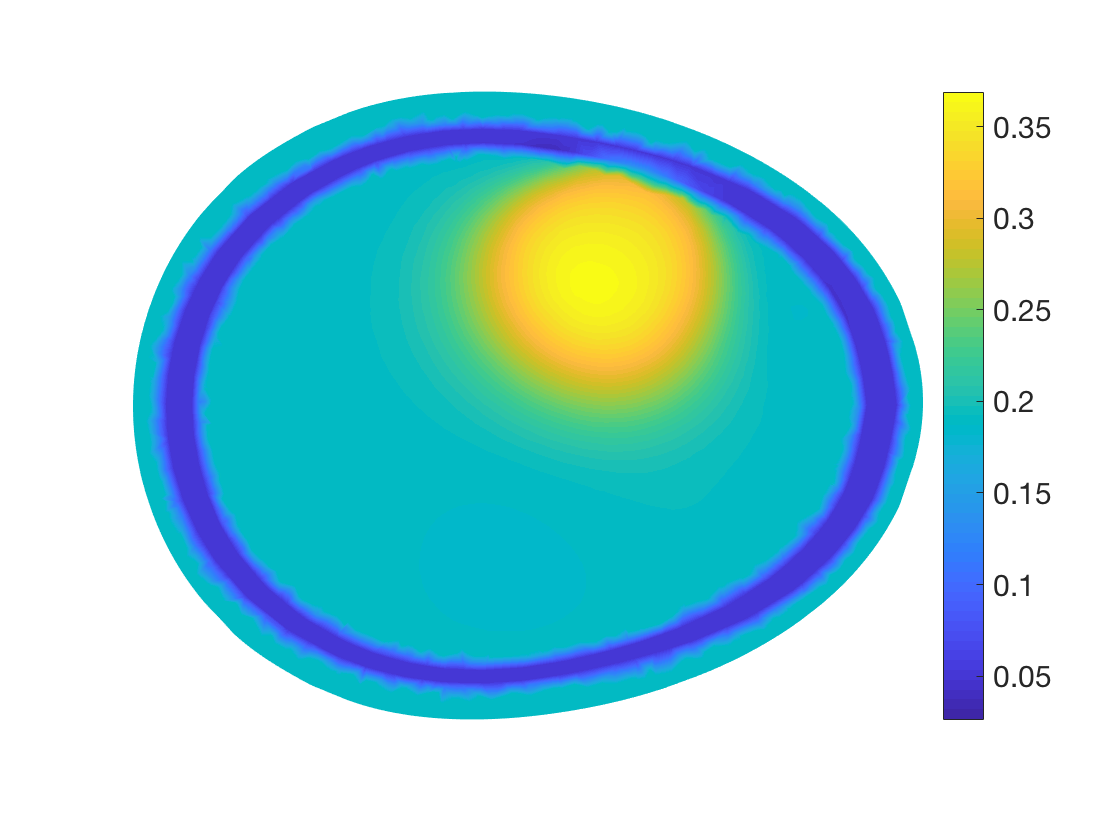}}
  {\includegraphics[width=4.1cm]{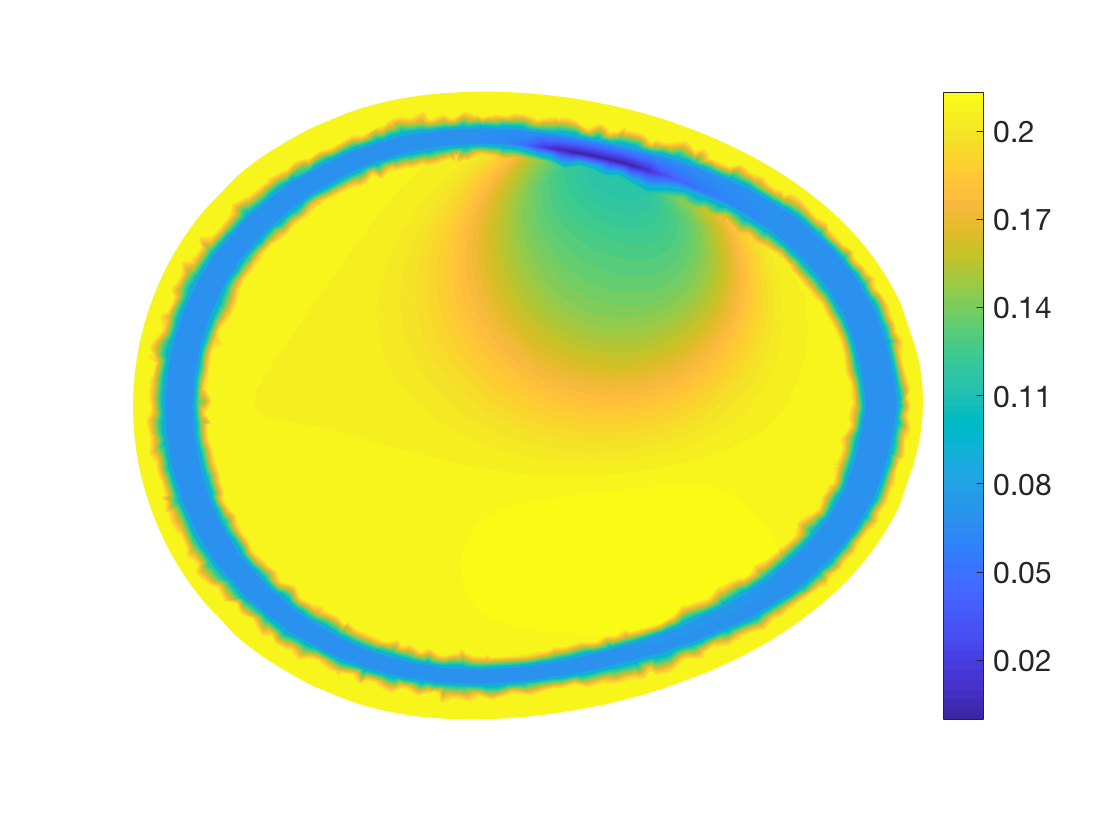}}
  }
    \caption{{\sc Case~1:} Slices at height $3$\,cm. Top row: target conductivities. Middle row: reconstructions without approximation error modeling. Bottom row: reconstruction with approximation error modeling. Left column: (i) healthy patient. Center column: (ii) hemorrhagic stroke. Right column: (iii) ischemic stroke.
    }
    \label{fig:case1}
\end{figure}

The first experiment tackles a setup where the (shape) parameters defining the target heads are relatively close to their average values; 
the results of the experiment are visualized in Figure~\ref{fig:case1}. The columns in Figure~\ref{fig:case1} correspond to the considered three medical conditions: a healthy patient on the left, a hemorrhagic stroke in the middle and an ischemic stroke on the right. The top row in Figure~\ref{fig:case1} shows cross sections of the conductivity inside the target heads at height $3$\,cm; recall that the hemorrhagic and the ischemic stroke manifest themselves as ball-shaped inclusions of radius $2.25$\,cm centered at $(2,3,3)$\,cm with conductivities $2$\,S/m and $0.02$\,S/m, respectively. The middle row presents slices from reconstructions produced by Algorithm~\ref{alg:kokohoska} without taking the approximation error noise into account, whereas the bottom row illustrates cross sections of reconstructions that are based on approximation error modeling. 

Notice that the colormaps in different subimages of Figure~\ref{fig:case1} are different: they have been chosen to cover the range of the corresponding (reconstructed) conductivity in the {\em whole} head. In particular, if the shown cross section of a reconstruction does not cover the entire color scale, the reconstruction in question contains artifacts not intersecting with the illustrated slice. This same conclusion applies also to reconstructions shown in Figures~\ref{fig:case2} and \ref{fig:case3} below. 

The reconstructions in the bottom row of Figure~\ref{fig:case1} do not suffer from significant artifacts, and the two on the right also clearly indicate the locations of the hemorrhagic and ischemic strokes. However, the conductivity levels of the inclusions are not reproduced accurately. These reconstructions also seem to be in a relatively good agreement with the TV prior: the reconstructed quantity, i.e.~the conductivity perturbation~$\kappa$, exhibits a relatively well localized inhomogeneity (if any) in an almost constant background. The dynamic range of the inclusions is arguably captured somewhat better when the approximation error noise is ignored in the middle row of Figure~\ref{fig:case1}, but the corresponding reconstructions suffer from considerable artifacts, and it would be difficult to infer what kinds of embedded inclusions there actually are inside the target head. In particular, the reconstruction corresponding to the healthy patient contains lots of small artifacts close to the interior surface of the skull.

Both with and without approximation error modeling the reconstruction of the hemorrhagic, i.e.~conductive, stroke is better than that of the ischemic, i.e.~insulating, one. This is a recurring theme throughout our numerical examples, and it is presumably caused by the difficulty in differentiating the impact of an ischemic stroke on the measurements from the shielding effect of the skull. In fact, assuming approximation error modeling is employed, one would be able to obtain reasonable reconstructions of even considerably smaller hemorrhagic strokes in all our numerical experiments, but one cannot, unfortunately, make a similar claim about ischemic strokes. However, the poorer resolution power of the reconstructions for ischemic strokes should not be a critical problem for the envisioned application of EIT stroke classification in emergency care, since it should be sufficient to be able to exclude a hemorrhage based on EIT measurements for a decision on whether to start clot dissolving medication for a patient with an acute stroke.

\subsubsection{Case~2:~head shape and other parameters further from the average}
\begin{figure}[b!]
\center{
  {\includegraphics[width=4.1cm]{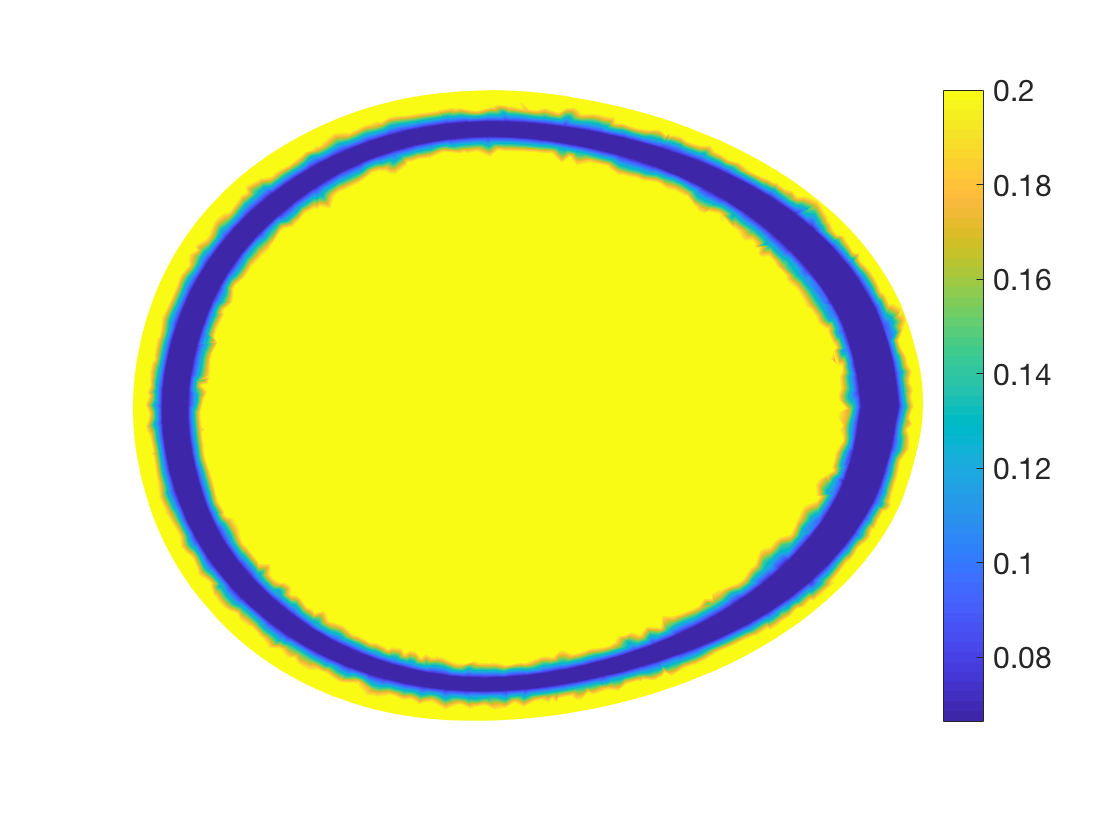}}
  {\includegraphics[width=4.1cm]{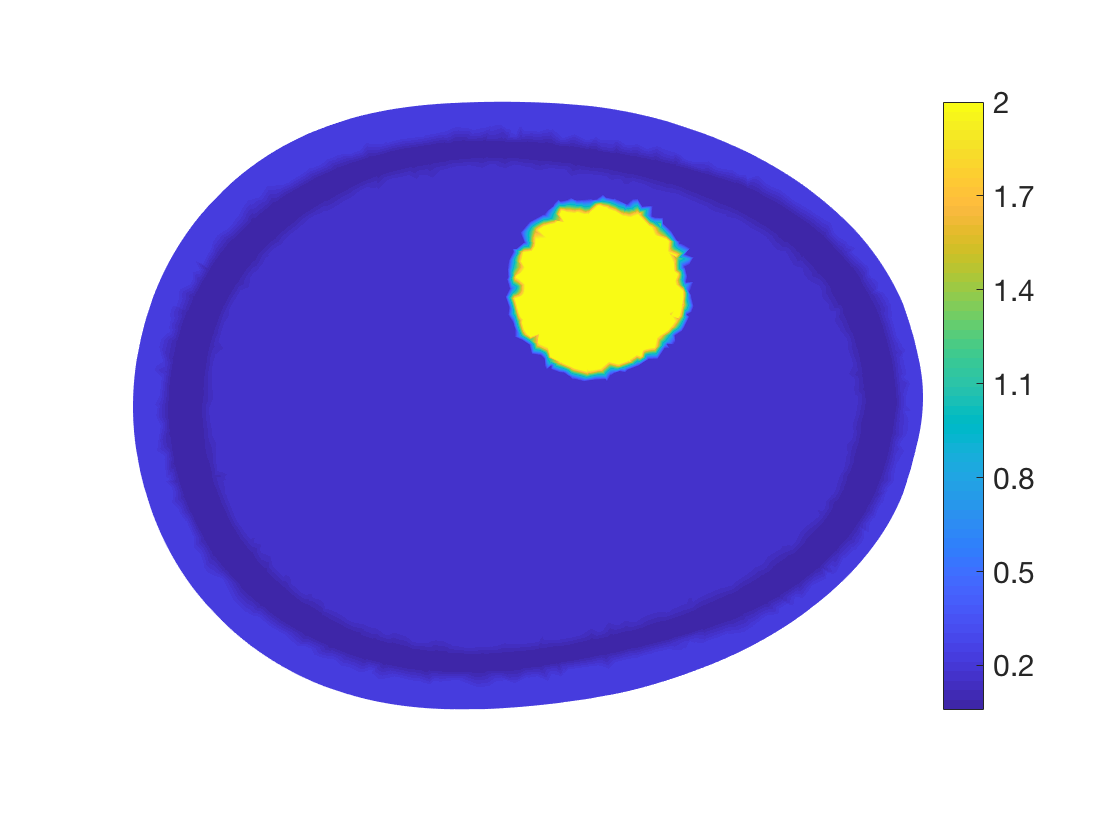}}
  {\includegraphics[width=4.1cm]{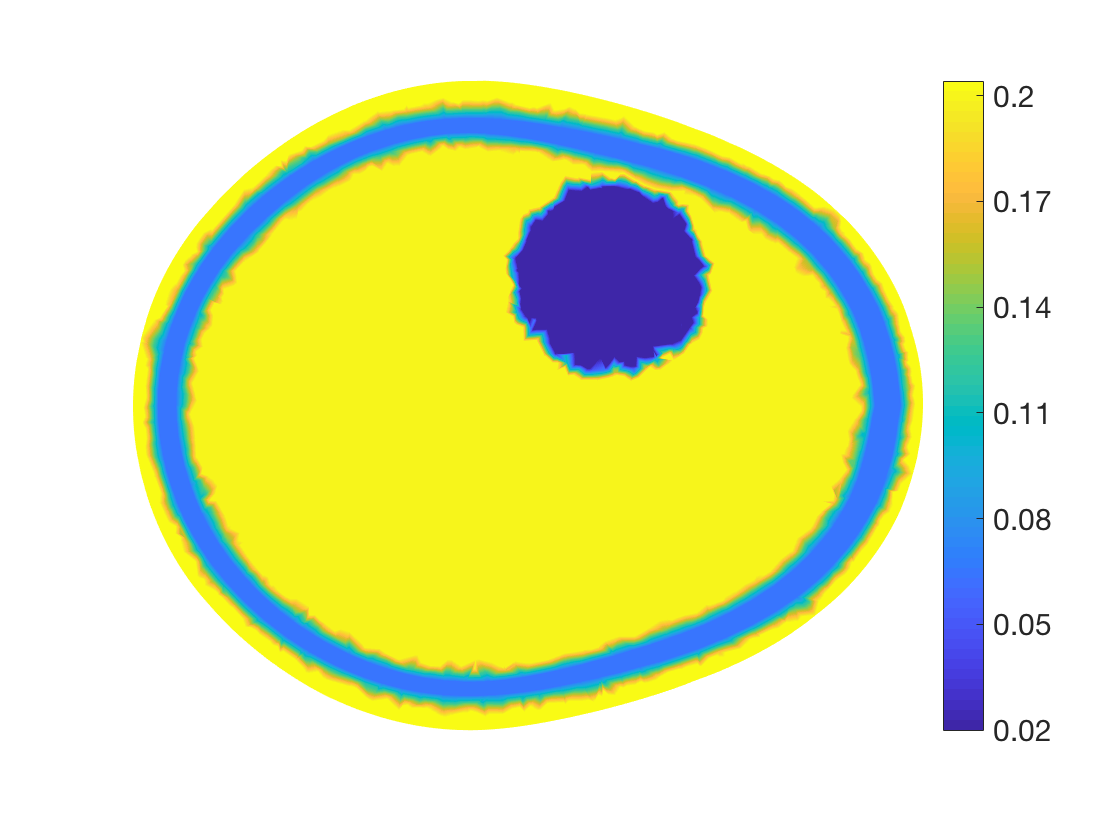}}
  }
  \center{
  {\includegraphics[width=4.1cm]{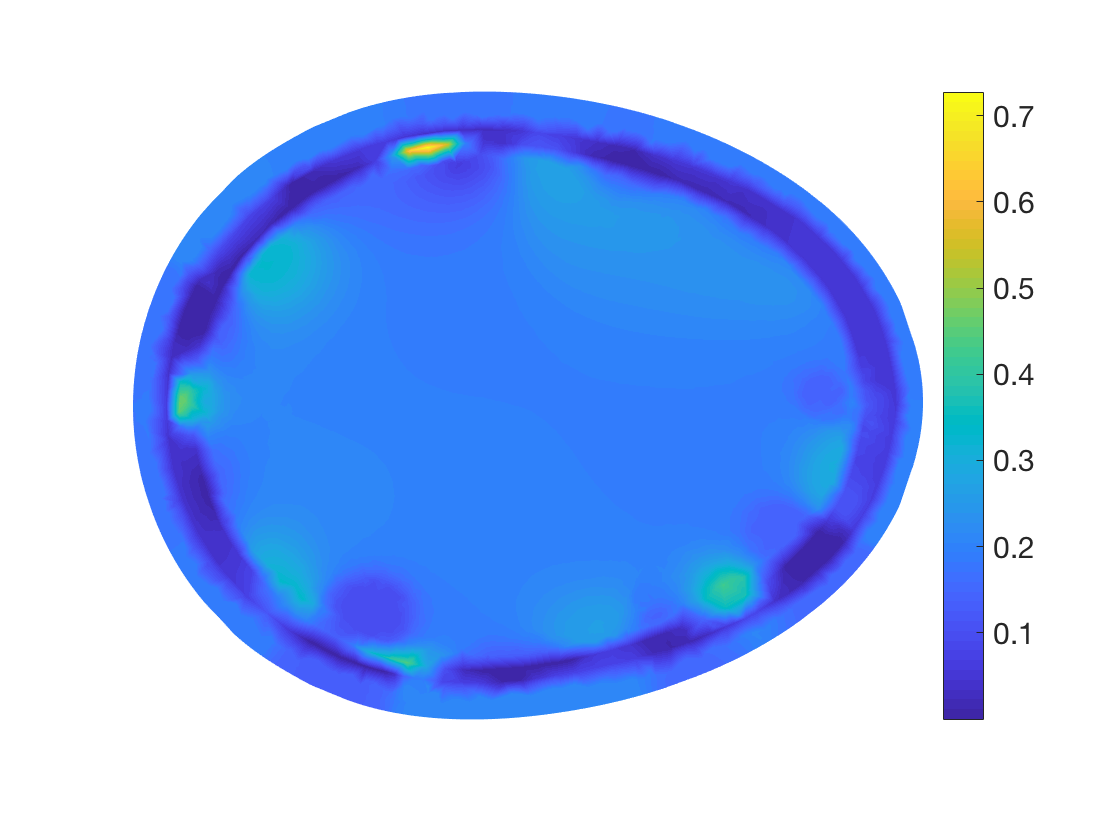}}
  {\includegraphics[width=4.1cm]{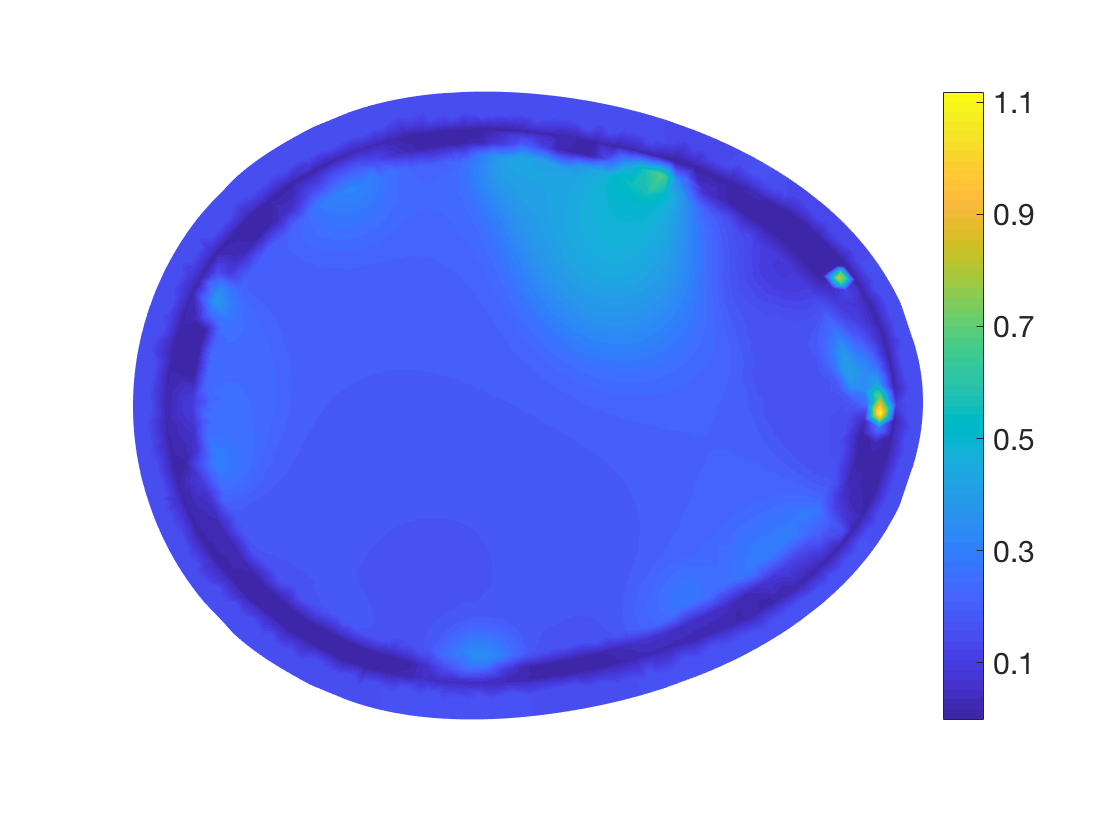}}
  {\includegraphics[width=4.1cm]{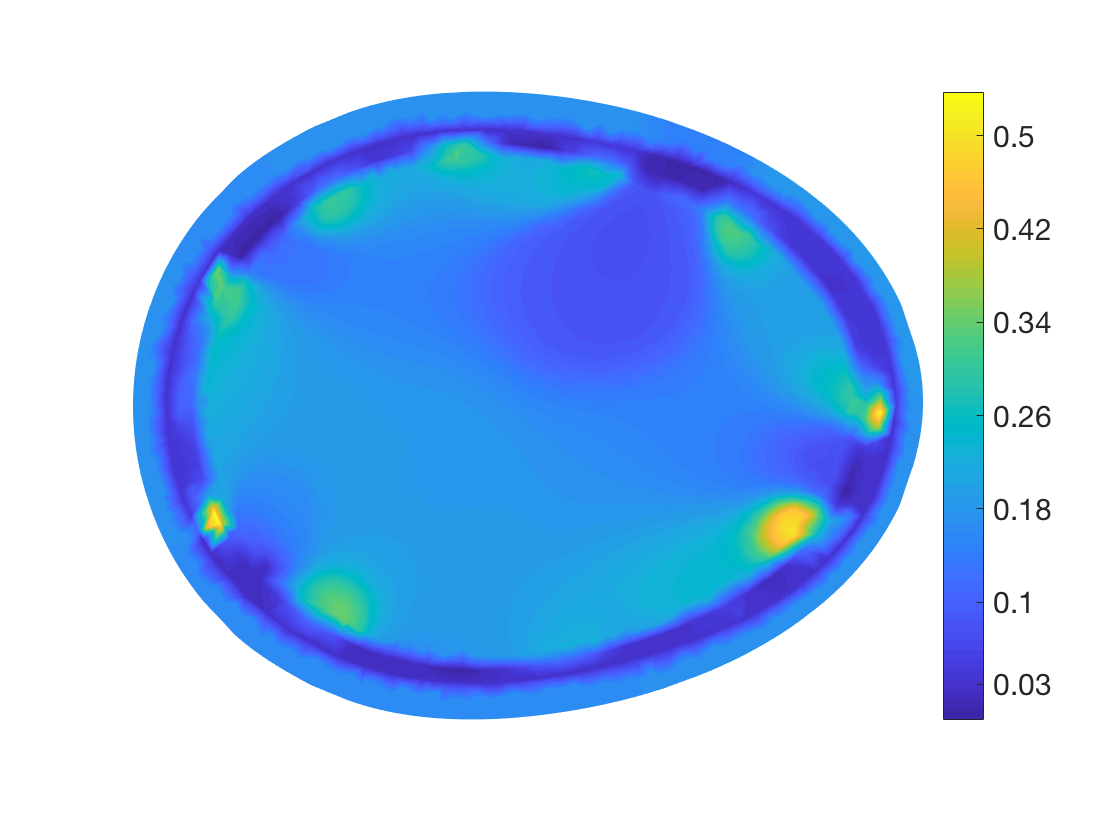}}
  }
  \center{
  {\includegraphics[width=4.1cm]{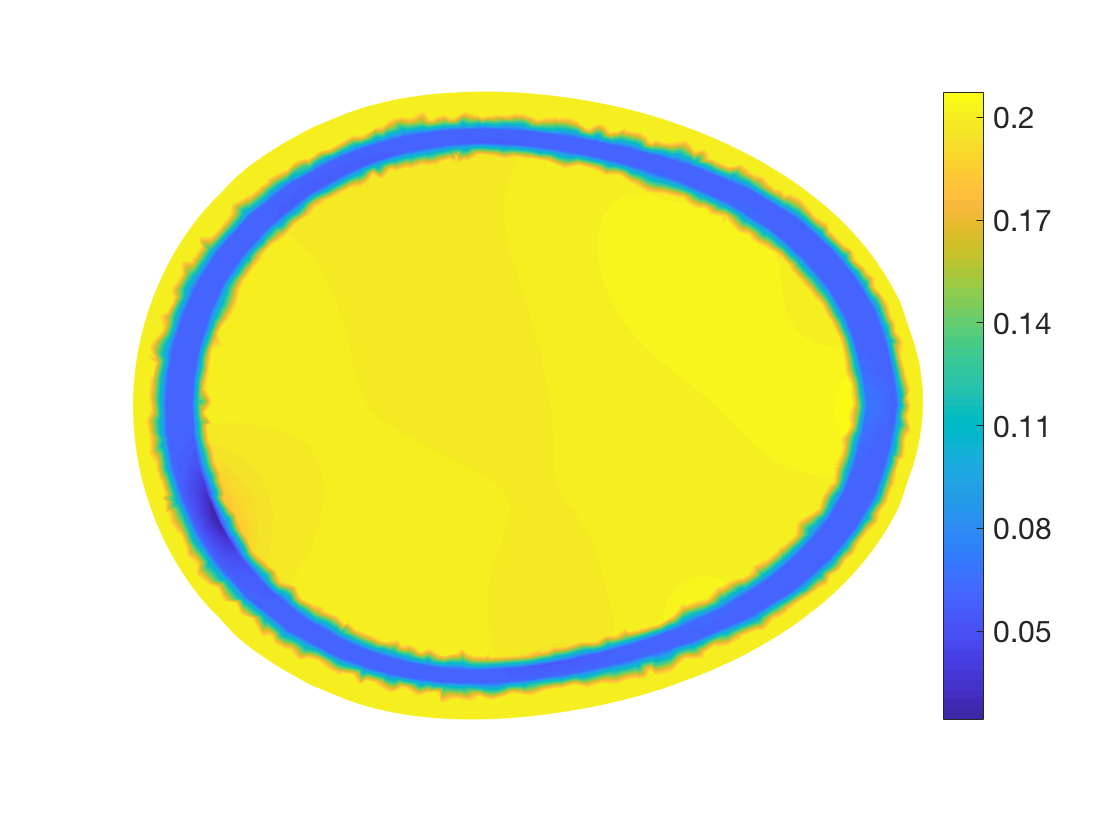}}
  {\includegraphics[width=4.1cm]{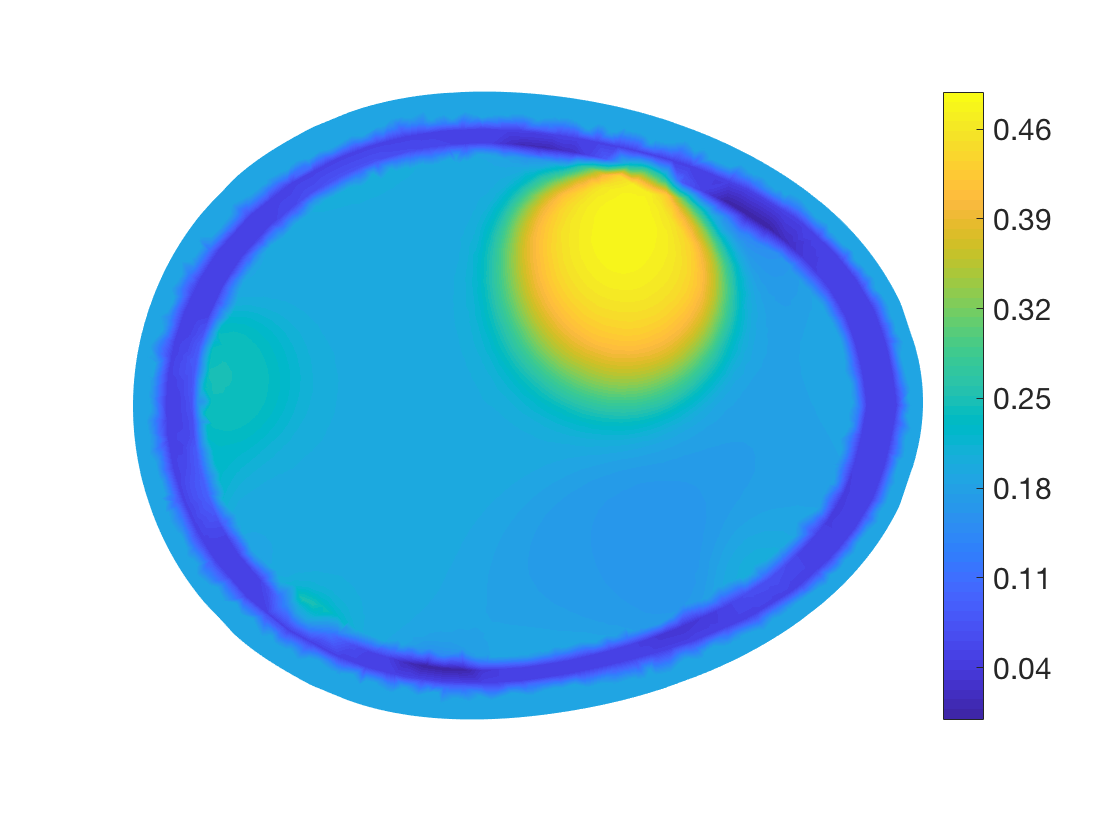}}
  {\includegraphics[width=4.1cm]{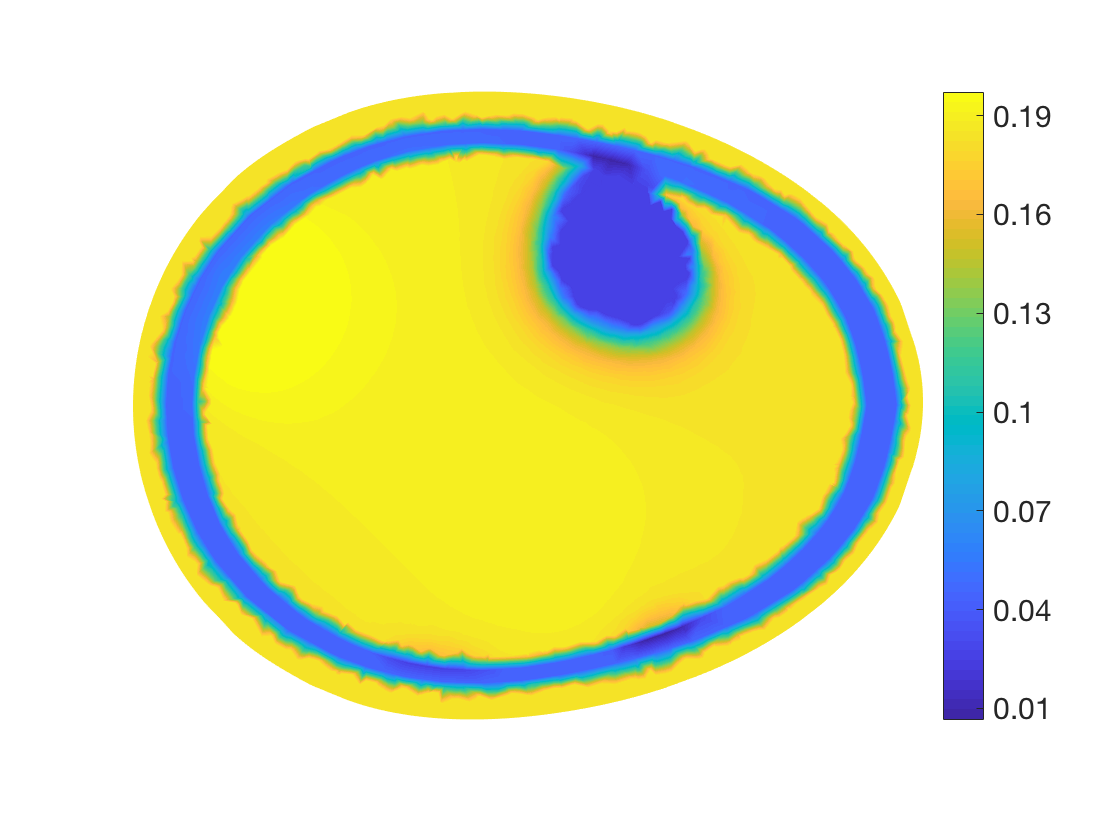}}
  }
    \caption{{\sc Case~2:} Cross sections at height $3$\,cm. Top row: target conductivities. Middle row: reconstructions without approximation error modeling. Bottom row: reconstruction with approximation error modeling. Left column: (i) healthy patient. Center column: (ii) hemorrhagic stroke. Right column: (iii) ischemic stroke.
    }
    \label{fig:case2}
\end{figure}

Our second experiment considers a setting where the parameters defining the targets are drawn from the distributions listed in Table~\ref{tab:values}; the results of the experiment are visualized in Figure~\ref{fig:case2}, which is organized in the same way as Figure~\ref{fig:case1} in Case~1. That is, the columns of Figure~\ref{fig:case2} correspond to the considered three medical conditions and its rows to the target head, reconstructions without approximation error modeling and reconstructions accounting for the approximation error noise, respectively.

This time the cross sections of the reconstructions without approximation error modeling in the middle row of Figure~\ref{fig:case2} are almost useless due to the severe artifacts located mainly close to the interior surface of the skull. On the other hand, the ones in the bottom row, accounting for the approximation error noise, are almost as informative as in Case~1. In particular, one can still identify and approximately locate the strokes in the two right-most images in the bottom row of Figure~\ref{fig:case2}, which is the main goal of the considered imaging application.

\subsubsection{Case~3:~average head and known parameters}
In our third and final test, we assume that the imaged patient has an average three-layer head and the contact resistances also take their expected values; in other words, the measurement data are simulated using the expected parameter values listed in Table~\ref{tab:values} and including an ischemic or hemorrhagic inclusion in the conductivity parametrization when appropriate. In particular, this choice makes the use of approximation error modeling unnecessary since the surrogate model and the accurate model are the same for the examined `average patient'.

\begin{figure}[b!]
\center{
  {\includegraphics[width=4.1cm]{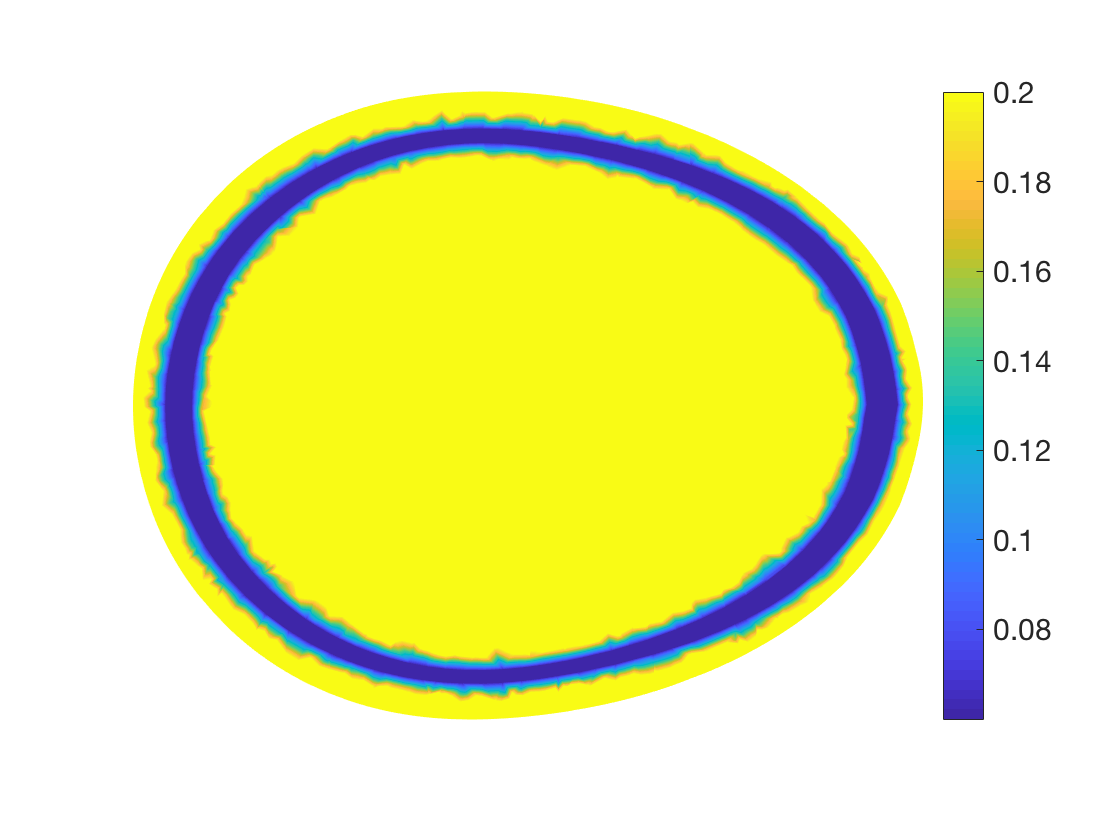}}
  {\includegraphics[width=4.1cm]{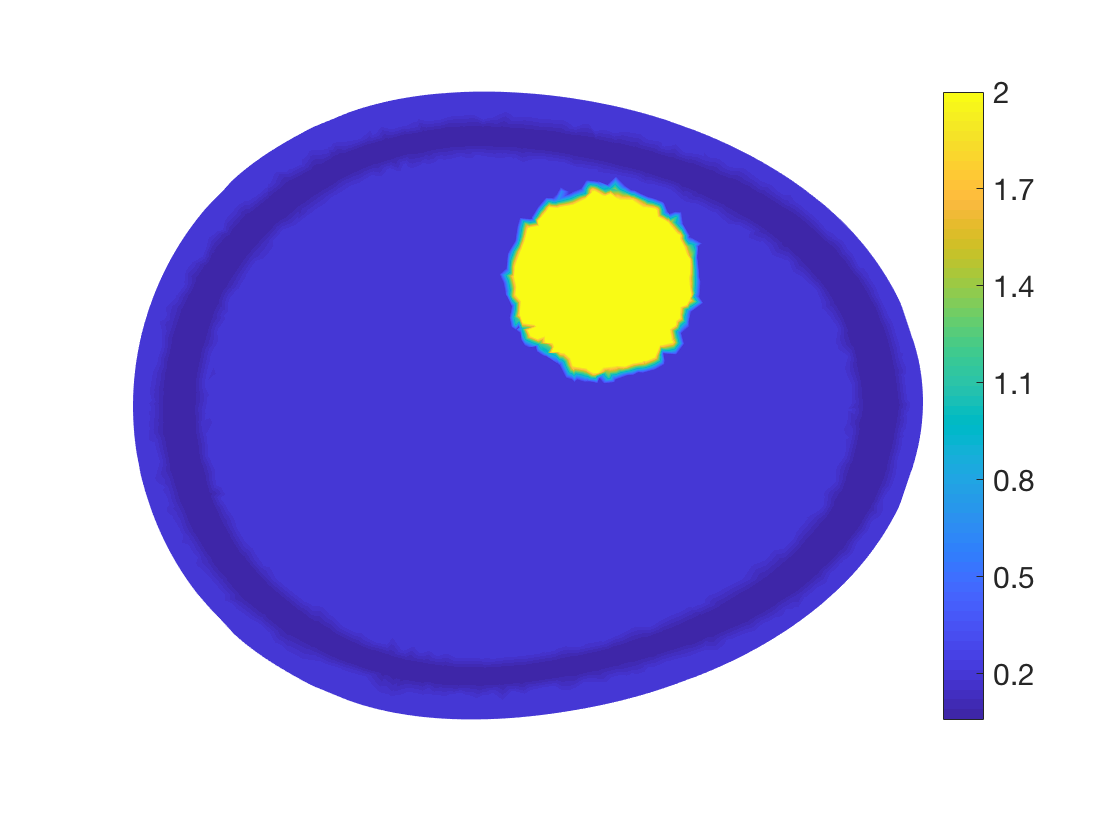}}
  {\includegraphics[width=4.1cm]{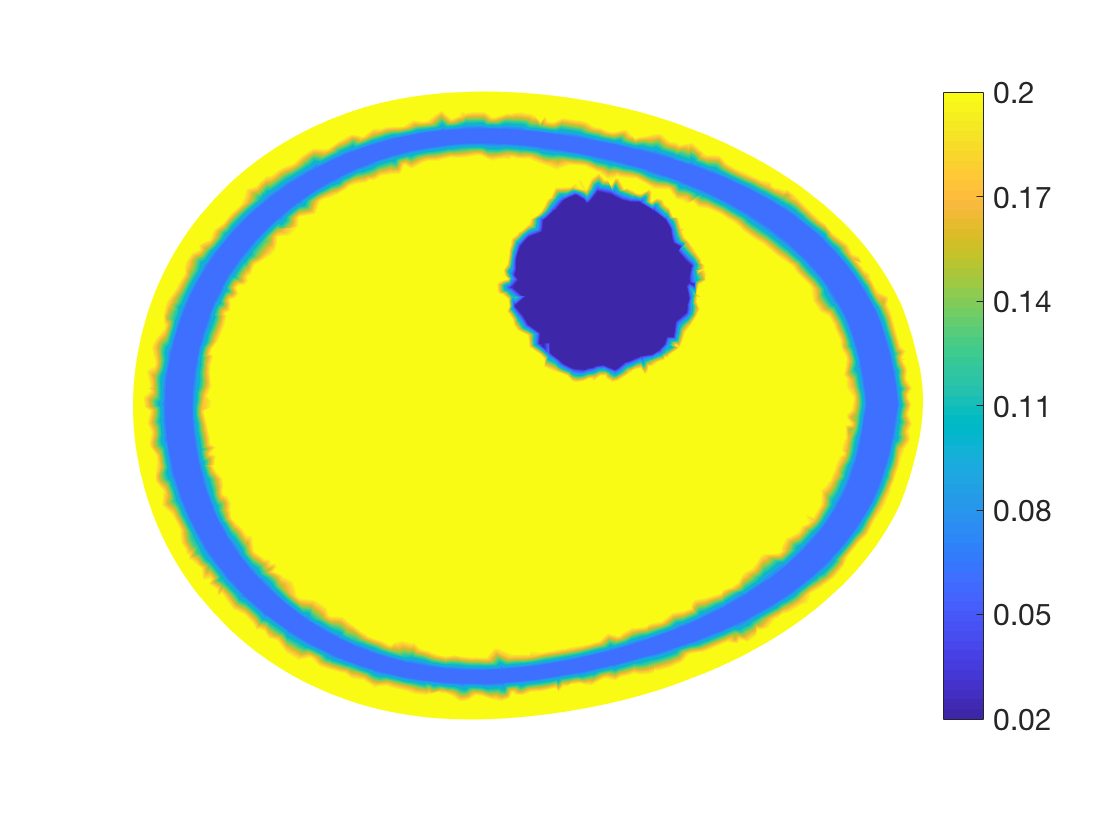}}
  }
  \center{
  {\includegraphics[width=4.1cm]{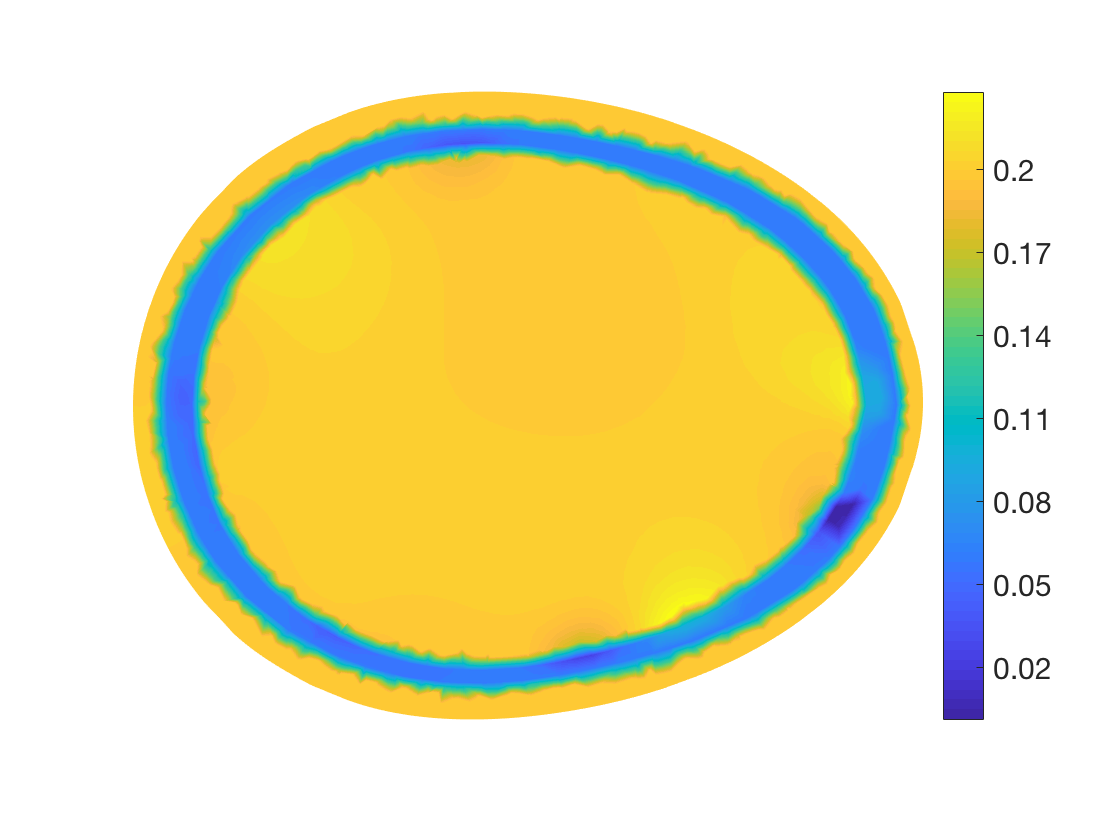}}
  {\includegraphics[width=4.1cm]{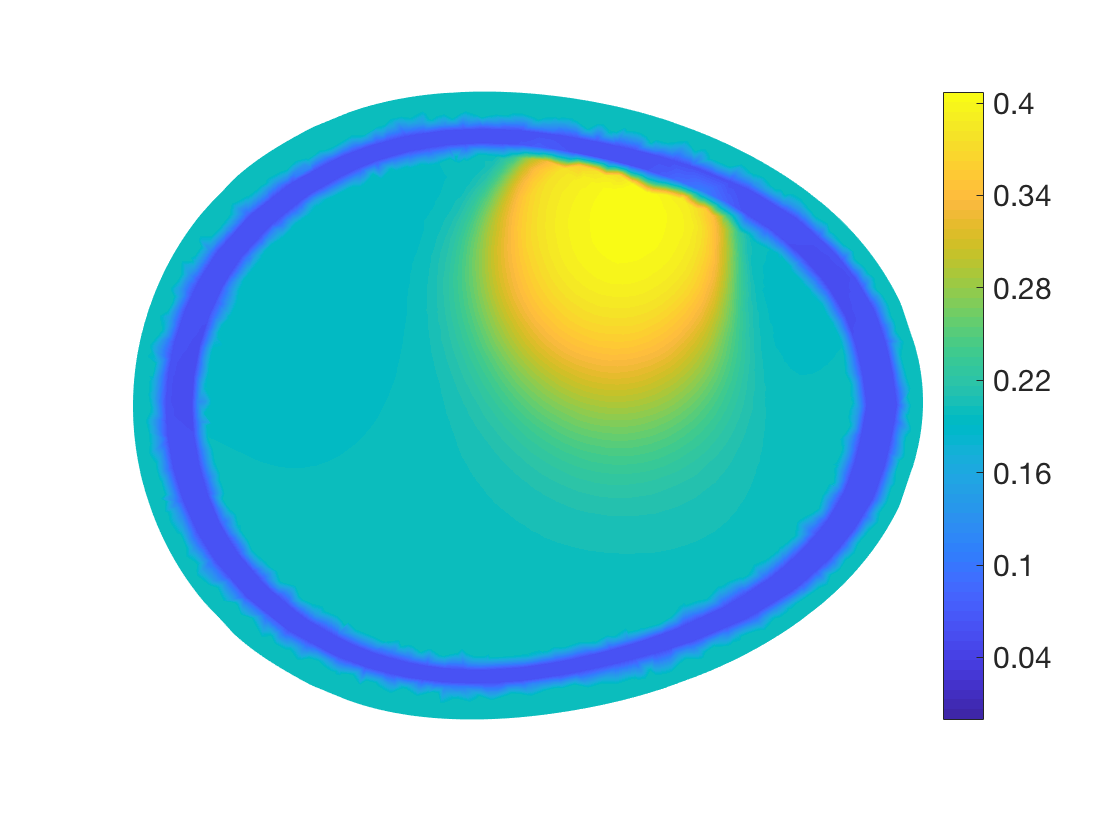}}
  {\includegraphics[width=4.1cm]{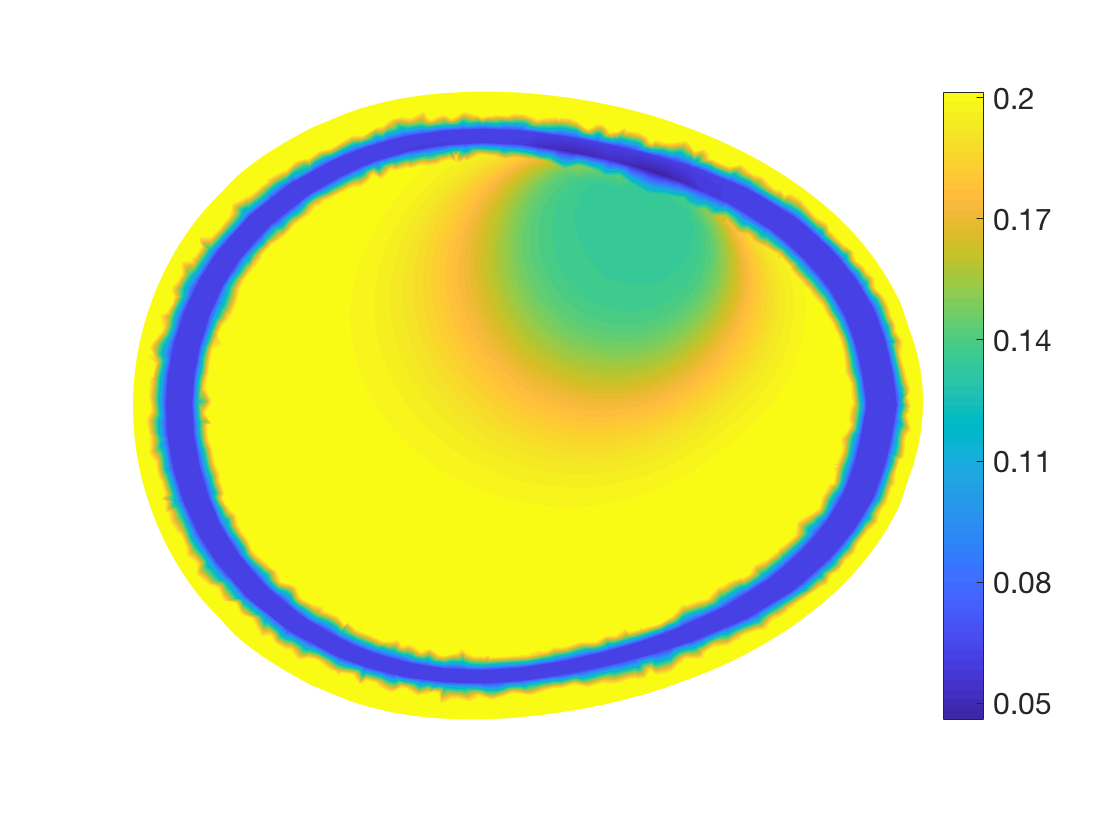}}
  }
  \center{
  {\includegraphics[width=4.1cm]{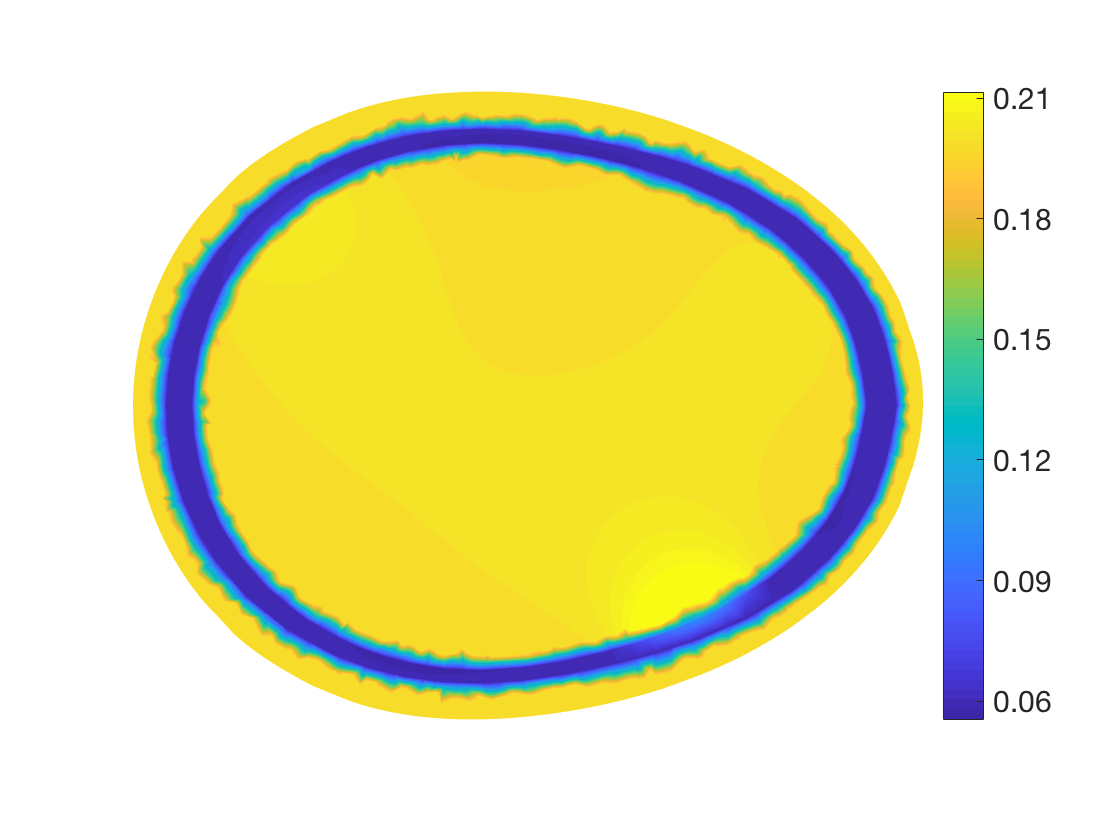}}
  {\includegraphics[width=4.1cm]{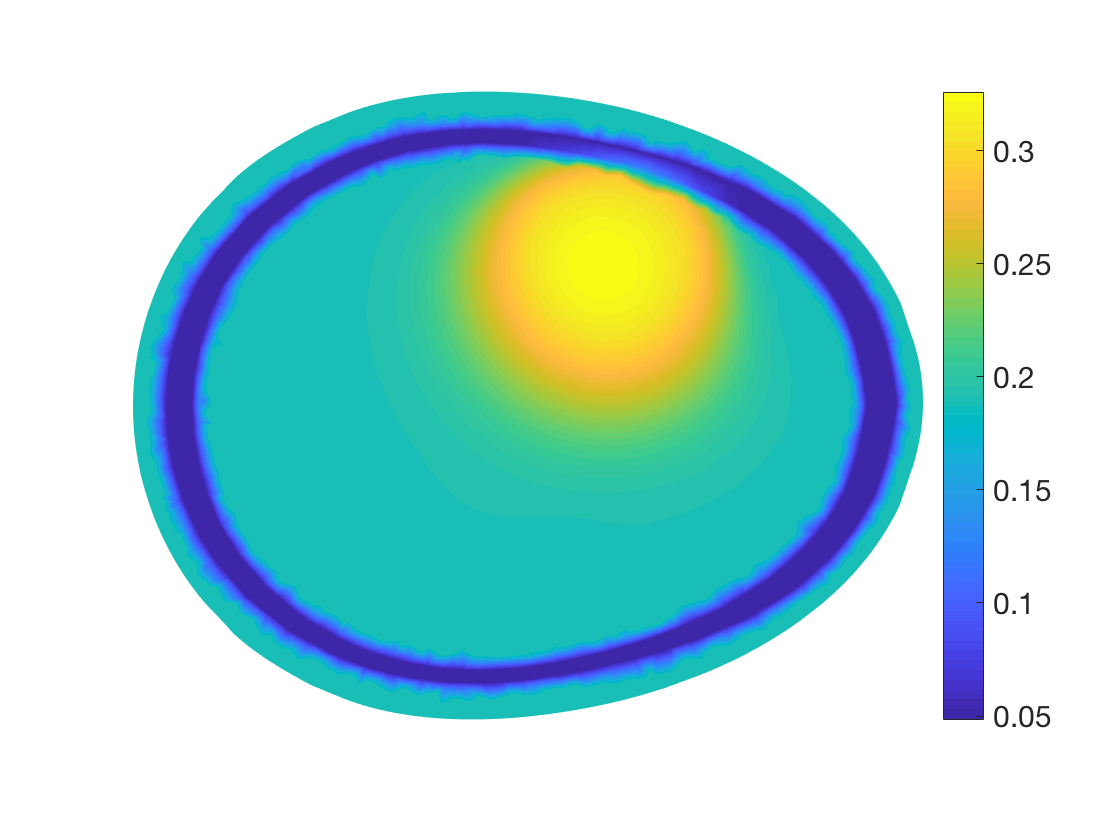}}
  {\includegraphics[width=4.1cm]{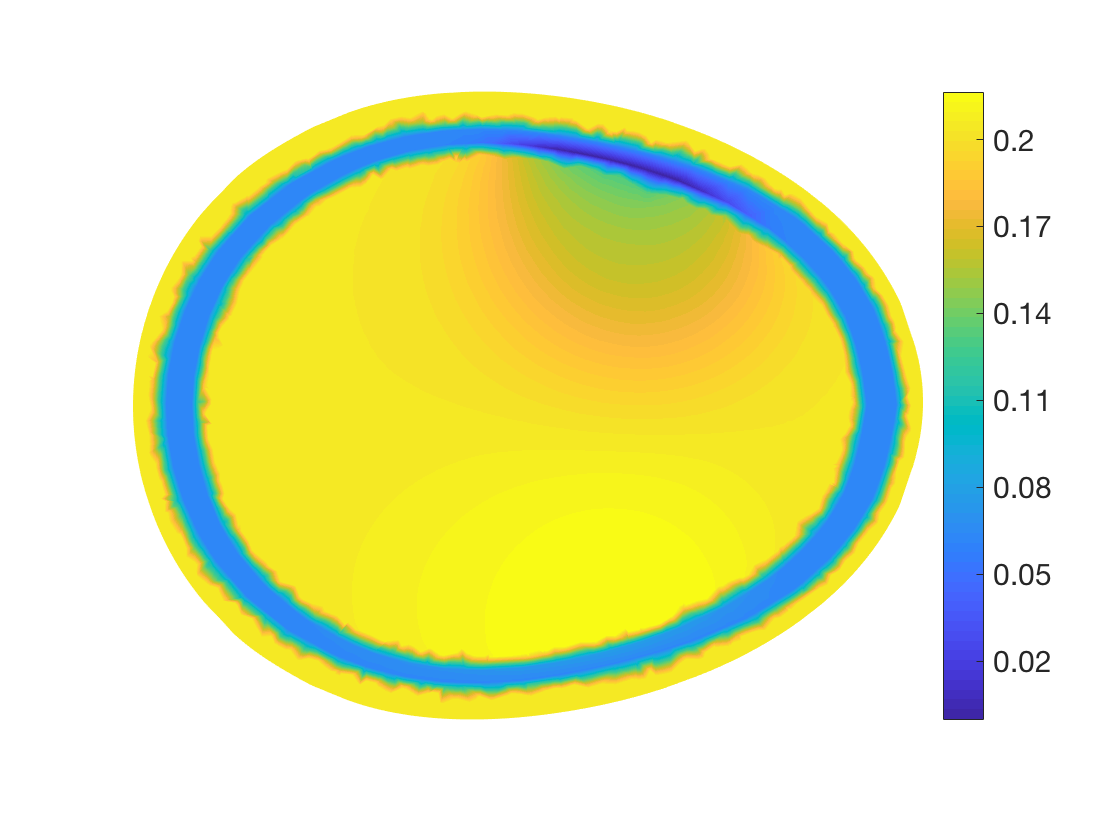}}
  }
    \caption{{\sc Case~3:} Cross sections at height $3$\,cm. Top row: target conductivities. Middle row: reconstructions without approximation error modeling. Bottom row: reconstruction with approximation error modeling. Left column: (i) healthy patient. Center column: (ii) hemorrhagic stroke. Right column: (iii) ischemic stroke.
    }
    \label{fig:case3}
\end{figure}

The results of the experiment are illustrated in Figure~\ref{fig:case3}, which is organized in the same way as the corresponding figures in the previous two cases. As expected, the reconstructions ignoring the approximation error noise in the middle row of Figure~\ref{fig:case3} are this time around somewhat better than those (unnecessarily) accounting for the (nonexisting) approximation error noise, with the difference being most notable for the ischemic stroke. However, all four reconstructions involving an embedded inclusion are well aligned with the TV prior, and all of them also reproduce the location of the stroke relatively accurately. It thus seems that adopting approximation error modeling does not considerably hinder the reconstruction quality even if there actually is no mismodeling in the measurement setup.

\section{Concluding remarks}
\label{sec:conclusion}

We have introduced a numerical algorithm for EIT based stroke classification in the presence of modeling errors in the electrode positions and the patient's head geometry. The method is based on the approximation error approach along with an edge-preferring prior and the computation of a MAP-like estimate following the ideas in \cite{Harhanen15}. Our imaging model corresponds to an average head of an anatomical atlas, and the approximation error model for the errors caused by the inaccurate knowledge of the electrode locations and the shapes of the head and the skull is constructed based on a probabilistic model induced by the atlas.

The method was tested with simulated electrode measurements from three-dimensional head models, with the highly resistive skull included, and compared to reconstructions with a conventional measurement error model. This comparison demonstrated that accounting for the approximation error model produces superior reconstructions compared to the conventional measurement noise model when errors in the electrode locations and the head geometry are present. Moreover, including the approximation error in the measurement model did not significantly hamper the reconstruction quality in the (unlikely emergency care) scenario where the imaging geometry is exactly known. 

Our results show that reconstructing the conductivity perturbation with respect to the prior mean in the average head model produces useful information on the stroke in the true measurement geometry, if the approximation error is taken properly into account.
These findings indicate that the proposed approach could be a feasible choice for stroke classification in emergency care --- an application where one would practically always lack the exact knowledge on the electrode locations and the geometry of the patient's head.

Our current three-layer head model has obvious limitations nonetheless, being a simplified version of the true head anatomy. In particular, it does not take into account the characteristic porosity of the skull or the {\em cerebrospinal fluid} (CSF), a layer of conductive liquid between the skull and the brain. The CSF has been previously considered within the approximation error paradigm in \cite{nissinen15}. Furthermore, the currently employed anatomical atlas covers (roughly) only the top half of the head. In real-life applications, the atlas may need to be replaced by a model that covers a larger portion of the head, or, alternatively, the domain truncation effects may need to be included by utilizing, e.g., stochastic boundary models~\cite{calvetti15a,calvetti15b,Hadwin14}. Moreover, we modeled hemorrhagic and ischemic strokes by high and low conductivity inclusions, with increased local conductivity corresponding to an hemorrhage (excess of highly conductive blood) and decreased conductivity to an ischemic stroke (lack of highly conductive blood). In real-life, the situation may be somewhat more complicated due to the presence of a penumbra, i.e.~a region of normal to high blood volume surrounding an ischemic stroke as the brain tries to balance the net blood pressure and flow, or hypodense tissue around a hemorrhage suffering from a shortage of blood. Such considerations and testing our algorithm with real-world stroke classification data are left for future studies.

\bibliographystyle{acm}
\bibliography{headae-refs.bib}

\end{document}